\documentclass[a4paper, 12pt, oneside, reqno, notitlepage]{amsart}
\usepackage{amsmath,amssymb,amsthm,graphicx,mathrsfs,bbm,url}
\usepackage[margin=3cm]{geometry}
\usepackage{amsthm}
\usepackage{wrapfig}
\usepackage{enumitem}
\usepackage{mathtools}
\usepackage[utf8]{inputenc}
 \usepackage[T1]{fontenc}
\usepackage[usenames,dvipsnames]{color}
\usepackage[colorlinks=true,linkcolor=Red,citecolor=Green]{hyperref}
\usepackage[super]{nth}
\usepackage[open, openlevel=2, depth=3, atend]{bookmark}
\hypersetup{pdfstartview=XYZ}
\usepackage[font=footnotesize]{caption}
\usepackage[colorinlistoftodos]{todonotes}
\usepackage[normalem]{ulem}

\usepackage{epstopdf}
 
\usepackage{hyperref}

\theoremstyle{plain}
\newtheorem{theorem}{Theorem}[section]
\newtheorem*{theorem*}{Theorem}

\newtheorem{lemma}[theorem]{Lemma}
\newtheorem{proposition}[theorem]{Proposition}
\newtheorem{corollary}[theorem]{Corollary}

\theoremstyle{definition}
\newtheorem{definition}[theorem]{Definition}
\newtheorem{example}[theorem]{Example}

\theoremstyle{remark}
\newtheorem{remark}[theorem]{Remark}

\numberwithin{equation}{section}

\newcommand{\C}{\mathbb{C}}
\newcommand{\R}{\mathbb{R}}

\newcommand{\Z}{\mathbb{Z}}
\newcommand{\D}{\mathbb{D}}

\newcommand{\M}{\mathcal{M}}
\newcommand{\N}{\mathbb{N}}
\newcommand{\T}{\mathbb{T}}
\newcommand{\V}{\mathbb{V}}

\newcommand{\HH}{\mathbb{H}}
\newcommand{\Ss}{\mathbb{S}}
\newcommand{\eps}{\varepsilon}

\newcommand{\mc}{\mathcal}

\newcommand{\dd}{\mathrm{d}}
\newcommand{\pol}{\mathrm{pol}}
\newcommand{\fib}{\mathsf{p}}

\DeclareMathOperator{\Op}{Op}
\DeclareMathOperator{\WF}{WF}

\DeclareMathOperator{\supp}{supp}

\DeclareMathOperator{\comp}{comp}
\DeclareMathOperator{\spn}{span}

\newcommand{\be}{\begin{equation}}
\newcommand{\ee}{\end{equation}}

% ...

% ...

\title
[Invariant distributions \& transport twistor spaces]
{Invariant distributions and the transport twistor space of closed surfaces}

\author[Bohr]{Jan Bohr}
\address{Mathematisches Institut
der Universit{\"a}t Bonn, Endenicher Allee 60, D-53115 Bonn }
\email{bohr@math.uni-bonn.de}

\author[Lefeuvre]{Thibault Lefeuvre}
\address{Université de Paris and Sorbonne Université, CNRS, IMJ-PRG, F-75006 Paris, France.}
\email{tlefeuvre@imj-prg.fr}

\author[Paternain]{Gabriel P. Paternain}
\address{Department of Pure Mathematics and Mathematical Statistics, University of
Cambridge, Cambridge CB3 0WB, UK, and Department of Mathematics, University of Washington, Seattle, WA 98195, USA.}
\email{g.p.paternain@dpmms.cam.ac.uk}

\begin{document}

\begin{abstract}
We study transport equations on the unit tangent bundle of a closed oriented Riemannian surface and their links to the \emph{transport twistor space} of the surface (a complex surface naturally tailored to the geodesic vector field). We show that fibrewise holomorphic distributions invariant under the geodesic flow  -- which play an important role in tensor tomography on surfaces -- form a \emph{unital algebra}, that is, multiplication of such distributions is well-defined and continuous. We also exhibit a natural bijective correspondence between fibrewise holomorphic invariant distributions and genuine holomorphic functions on twistor space with polynomial blowup on the boundary of the twistor space. Additionally, when the surface is Anosov, we classify holomorphic line bundles over twistor space which are smooth up to the boundary. As a byproduct of our analysis, we obtain a quantitative version of a result of Flaminio \cite{Flaminio-92} asserting that invariant distributions of the geodesic flow of a positively-curved metric on $\Ss^2$ are determined by their zeroth and first Fourier modes.
\end{abstract}


\maketitle
\setcounter{tocdepth}{1}

%

%
%
%
%
%
%
%
%

\section{Introduction}


Let $(M,g)$ be a closed oriented Riemannian surface with unit tangent bundle $SM$ and geodesic vector field $X$, considered as a first-order differential operator $X\colon C^\infty(SM)\rightarrow C^\infty(SM)$. This paper is concerned with the  {\it transport equation}
\begin{equation}\label{TE1}
Xu  =f
\end{equation}
and its interplay with {\it vertical Fourier analysis}. The latter refers to the decomposition $f=\sum_{k\in \Z} f_k$ of a function $f\in C^\infty(SM)$ into its vertical Fourier modes.  This corresponds to freezing $x\in M$ and writing $f(x,\cdot)\colon S_x M\simeq \mathbb S^1\rightarrow \C$ as Fourier series in the velocity variable. Globally, the $k$th Fourier mode $f_k$ may be viewed either as a function in $C^\infty(SM)$ or as a section of a complex line-bundle $\Omega_k\rightarrow M$, see \textsection \ref{section:fourier}.

As observed by Guillemin and Kazdhan \cite{Guillemin-Kazhdan-80}, the transport equation  \eqref{TE1}
is  equivalent to a system of equations for the Fourier modes
\begin{equation}\label{TE2} \eta_+ u_{k-1} + \eta_- u_{k+1} = f_k,\quad k\in \Z,\end{equation} where  $\eta_\pm \colon  C^\infty(M,\Omega_k)\rightarrow C^\infty(M,
\Omega_{k\pm 1})$
 are elliptic differential operators. Starting with the just cited article, the interplay between transport equations and vertical Fourier analysis has become a central feature in the theory of geometric inverse problems in two dimensions \cite{Paternain-Salo-Uhlmann-book}. Indeed, many such problems can be rephrased as assertions about the Fourier support $\{k\in \Z:f_k\neq 0\}$ and frequently it is key to produce solutions to the transport equation \eqref{TE1} whose Fourier modes satisfy additional constraints. An important example is the constraint
  $f_k =0$ for $k<0$ in which case $f$ is called {\it fibrewise holomorphic} (see also Definition \ref{def:fibrehol} below).
 
 This is epitomised by the {\it tensor tomography problem} for Anosov surfaces.
  For tensors of degree $m\in \N$ this is equivalent to showing, for $u\in C^\infty(SM)$ and $f=Xu$,  that
 \begin{equation}
f_k=0 \text{ for } \vert k \vert 
\ge m+1 \quad \Longrightarrow \quad \begin{cases} 
u_0 = \text{const.} & \text{if } m=0,\\
u_k =0  \text{ for }\vert k \vert \ge m, & \text{if } m\ge 1.
\end{cases}
 \end{equation}
Using the Pestov energy identity, this was proved for $m=0$ and $m=1$  in \cite{Dairbekov-Sharafutdinov-03} and, assuming negative curvature, also for $m\geq 2$ in \cite{Guillemin-Kazhdan-80}. For general Anosov surfaces and $m= 2$, the tensor tomography problem was solved in \cite{Paternain-Salo-Uhlmann-14-2}, and for all $m\geq 2$ in \cite{Guillarmou-17-1}. A key insight of the latter two articles was that the problem is equivalent to the existence of sufficiently many fibrewise holomorphic (distributional) solutions to $Xu=0$. One purpose of this article is to revisit this aspect from the point of view of `transport twistor spaces', thereby shedding new light on the tensor tomography problem and  providing an alternative proof.

Twistor spaces associated to projective structures have been considered by several authors starting with \cite{Hit80,LeB80}, see also \cite{Dub83,OR85,LeMa02,LM10,MePa20,Met21}. One common description of this twistor space is as the bundle $\mathcal J$ over $M$ whose fibre at $x\in M$ consists of the orientation compatible linear complex structures on $T_{x}M$.
In this paper, we consider the {\it transport twistor space} as described in \cite{Bohr-Paternain-21}. This a 2--1 branched cover over (the compactification of) $\mathcal J$ 
and the just cited description 
is particularly tailored  to the geodesic vector field; further, it is suitable to capture {\it all} invariant distributions of interest in dynamics and geometric inverse problems, while $\mathcal J$ misses half of them.

Transport twistor spaces allow to study the interplay between transport equations and vertical Fourier analysis through the lens of complex geometry.
 For any oriented Riemannian surface $(M,g)$, possibly with a non-empty boundary, the transport twistor space 
is the $4$-manifold
\begin{equation}
Z=\{(x,v)\in TM:g(v,v)\le 1\},
\end{equation}
equipped with a natural complex structure that turns the interior $Z^\circ$ into a classical complex surface, but which degenerates at $SM\subset \partial Z$ in a way that encodes the transport equation \eqref{TE1}---we refer to \textsection\ref{section:twistor} for precise definitions. In \cite{Bohr-Paternain-21} several twistor correspondences were set up that relate the transport equation 
to complex geometric objects. For example, the algebra of holomorphic functions that are smooth up to the boundary $\partial Z$, 
  \begin{equation}\label{smoothalgebra}
\mathcal A(Z)=\{f\in C^\infty(Z): f \text{ holomorphic in } Z^\circ\},
\end{equation}
was shown to be isomorphic to the space of smooth fibrewise holomorphic first integrals of the geodesic flow; this is implemented by the following map:
\begin{equation}\label{smoothcorrespondence}
	\mathcal A(Z) \xrightarrow{\sim} \{u\in C^\infty(SM): Xu = 0, ~u \text{ is fibrewise holomorphic} \},\quad f\mapsto f\vert_{SM}.
\end{equation}
While the mentioned correspondence principles hold true on general oriented Riemannian surfaces (with or without boundary), they are of limited interest in the closed case. 
 Indeed, if the geodesic flow of $(M,g)$ is {\it ergodic}, then there are no non-constant smooth (or even $L^2$-integrable) solutions to $Xu=0$, such that  $\mathcal A(Z)\cong \C$.
 
The focus of this article lies on the twistor space of closed surfaces. We first extend the correspondence principle \eqref{smoothcorrespondence}
 to holomorphic functions that are allowed to blow up at the boundary $\partial Z$. This allows to also capture distributional solutions to the transport equation, which play an important role in the context of the  tensor tomography problem. We then go on to classify
 holomorphic line bundles on $Z$, considered up to isomorphism on $Z^\circ$, by exhibiting a finite-dimensional moduli space---this complements the main result of \cite{Bohr-Paternain-21} which states that the twistor space of a simple surface (a simply connected surface with strictly convex boundary and no conjugate points) supports no non-trivial holomorphic vector bundles.
 

At this point, the main purpose of transport twistor spaces is of conceptual nature. For example, while the main theorem of \cite{Bohr-Paternain-21} is equivalent
to an Oka--Grauert type result for the twistor space $Z$ and while it was inspired by this twistorial point of view, its proof relied in an essential way on the theory of transport equations and microlocal analysis---the same holds true for the results of this article. This is in contrast with the main result in \cite{LeMa02} where the full force of a geometric result is used via a desingularisation of the twistor space $\mathcal J$ (namely, that $\mathbb{CP}^{2}$ has a unique complex structure up to biholomorphism). While it remains an intriguing question whether  tools from complex geometry can be used also in the present setting, we believe that already exposing this connection with complex geometry provides valuable additional insights, as exemplified by Theorem \ref{mainthm1} below. Moreover, it is worthy of note that microlocal analysis and hyperbolic dynamics on the unit tangent bundle allow to describe complex geometric features of the transport twistor space as illustrated by Theorem \ref{modulianosov-intro} below.

\subsection{Algebras of holomorphic functions}\label{intro:algebras} Let $Z$ be the twistor space of a closed, oriented surface $(M,g)$. We consider the following three algebras of holomorphic functions:
\begin{equation}
 \mathcal A(Z) \subset \mathcal A_\pol(Z) \subset \mathcal A(Z^\circ)
\end{equation}
The largest algebra $\mathcal A(Z^\circ)$ contains all holomorphic functions $Z^\circ\rightarrow \C$, with no restriction on their boundary behaviour. Next, $\mathcal A_\pol(Z)$ contains all $f\in \mathcal A(Z^\circ)$ that satisfy a polynomial growth bound at $\partial Z$ as follows: there exist constants $C,p>0$ such that 
\begin{equation}
\label{equation:blowup-pol}
\sup_{(x,v)\in SM} \vert f(x,rv) \vert \le C (1-r)^{-p}\quad \text{ for all }0<r<1.
\end{equation} 
Finally, $\mathcal A(Z)=\mathcal A(Z^\circ)\cap C^\infty(Z)$ is the algebra from \eqref{smoothalgebra}.
Further, on the closed $3$-manifold $SM=\partial Z$ we consider the space of distributions $\mathcal D'(SM)$. We note that notions from vertical Fourier analysis such as fibrewise holomorphicity carry over from smooth functions to distributions.



\begin{theorem}[Twistor correspondence for invariant distributions] \label{mainthm1} There is a well-defined trace map $f\mapsto f\vert_{SM}$ from $\mathcal A_\pol(Z)$ into $\mathcal D'(SM)$, such that
\begin{equation}
\mathcal A_\pol(Z)\xrightarrow{\sim} \{u\in \mathcal D'(SM): Xu =0, ~u \text{ is fibrewise holomorphic}\},\quad f\mapsto f\vert_{SM}
\end{equation}
is an isomorphism. In particular, the space of fibrewise holomorphic invariant distributions forms a unital algebra. 
\end{theorem}

This result is reminiscent of the classical representation of distributions on $\R$ as boundary values of holomorphic functions on the upper halfplane in $\C$, see \cite[Theorems 3.1.11, 3.1.14]{Hormander-90}, where e.g.~the distribution $(x+i0^+)^{-1}$ can be viewed as distributional trace of $1/z$.

The fact that fibrewise holomorphic invariant distributions are closed under multiplication regardless of the underlying geometry is somewhat surprising. Indeed, multiplying such distributions is an important step in the solution of the tensor tomography problem for Anosov surfaces \cite{Paternain-Salo-Uhlmann-14-1, Guillarmou-17-1} and
so far it has been assumed that this was delicate issue. In the cited papers the multiplication problem was therefore only considered in particular instances, and their techniques heavily depended on the Anosov property. Using Theorem \ref{mainthm1} to justify multiplication thus allows to simplify the proofs in \cite{Paternain-Salo-Uhlmann-14-1, Guillarmou-17-1}, see also Remark 9.4 in the former paper; in particular, the more intricate analysis of Schwartz kernels in \cite{Guillarmou-17-1} can be circumvented.

If $u\in \mathcal D'(SM)$ is a fibrewise holomorphic  invariant distribution, then \eqref{TE2} together with the ellipticity of $\eta_-$ immediately implies that the Fourier modes of $u$ are all smooth. Any lack of regularity of $u$ is thus due to the growth of the norms $\Vert u_k\Vert_{C^N(M,\Omega_k)}$ as $k\rightarrow \infty$. Moreover, defining products of such distributions and extending $u$ to a function in $\mathcal A_\pol(Z)$ both hinge on controlling the growth of these norms. In order to do this, we analyse the wavefront set $\WF(u)$ and show that it is always contained in a  fixed convex cone inside $T^*(SM)\backslash 0$. This also shows directly that the multiplication 
is well defined by the classical multiplication theorem for distributions.

The analysis of the wavefront set of fibrewise holomorphic invariant distributions combines in a novel way the twist property of the geodesic flow with the microlocal propagation of singularities and the resulting regularity statement (see Theorem \ref{theorem:regularity} and Corollary \ref{corollary:freeconjugate} below) is of independent interest. It gives in particular that any fibrewise holomorphic invariant distribution as in Theorem \ref{mainthm1} has singular support free of conjugate points.

\medskip

For the next theorem, we introduce some additional notation: first, $\mathcal A(\D)\subset \mathcal A_\pol(\D)\subset\mathcal A(\D^\circ)$ are algebras of holomorphic functions on the compact unit disk $\D$, with boundary behaviour (smooth, polynomial blow-up or unrestricted) defined analogously to $Z$. Further, we turn $\mathcal A_\pol(Z)$ into a graded algebra by defining $\mathcal A^m_\pol(Z)$ $(m\ge 0)$ as the space of those functions $f\in \mathcal A_\pol(Z)$ whose trace $u=f\vert_{SM}$ satisfies $u_k=0$ for $k<m$.  


\begin{theorem}[Computation of some algebras of holomorphic functions]\label{mainthm2}
The following holds:
\begin{enumerate}[label=\emph{(\roman*)}]
\item If $M$ is diffeomorphic to $\mathbb S^2$, then
\[
\mathcal A(Z) = \mathcal A_\pol(Z) =\mathcal A(Z^\circ) \simeq \C.
\]
\item If $(M,g)$ is isometric to a flat torus, then
\begin{equation*}
\mathcal A(Z) \simeq \mathcal A(\D),\quad \mathcal A_\pol(Z) \simeq \mathcal A_\pol(\D),\quad \mathcal A(Z^\circ) \simeq \mathcal A(\D^\circ)
\end{equation*}
and these isomorphisms are induced by the inclusion $\D\hookrightarrow Z$ as a fibre.
\item If $(M,g)$ is Anosov, then 
\[
	\mathcal A(Z) \simeq \C,
\]
and both $\mc{A}_\pol(Z), \mc{A}(Z^\circ)$ are infinite-dimensional. Assuming moreover that the underlying Riemann surface $(M,[g])$ is non hyperelliptic,  the spaces $\mathcal A_\pol^m(Z)$ ($m\ge 0$) fit into a short exact sequence
\begin{equation*}
0\rightarrow \mathcal A_\pol^{m+1}(Z) \rightarrow \mathcal A_\pol^{m}(Z) \rightarrow \mathcal H_m \rightarrow 0,\quad m\ge 0,
\end{equation*}
where $\mathcal H_m$ is the space of holomorphic $m$-differentials of $(M,[g])$ (The condition on hyperellipticity can be dropped in the case $m=0$ and $m=1$). 

\end{enumerate}
\end{theorem}

It was suggested in \cite{Bohr-Paternain-21} that the twistor space of a simple surface behaves like a contractible Stein surface, and for a Euclidean domain $M\subset \R^2$ it was shown that the interior $Z^\circ$ is indeed Stein. For a closed surface $(M,g)$ on the other hand,  $Z^\circ$ can {\it never} be a Stein surface by the following elementary observation: the $0$-section embedding $\iota_M\colon M\rightarrow Z^\circ$ is holomorphic and hence the pull-back of any function in $\mathcal{A}(Z^\circ)$ is constant on $M$, which means that the algebra $\mathcal A(Z^\circ)$ cannot separate points on $M$.
Nevertheless, we may ask whether $Z^\circ$ is {\it holomorphically convex}. By the preceding theorem, the answer is {\it no}  if $M\approx \mathbb S^2$ and it is {\it yes} if $(M,g)$ equals a flat torus (see also Proposition \ref{proposition:t2}). In the Anosov case, where the algebra $\mathcal A_\pol(Z)$ appears to be fairly large, we leave this as an open question:

\medskip

\noindent{\bf Question.}~Does the twistor space of an Anosov surface have a holomorphically convex interior?



\subsection{Regularity of invariant distributions on convex surfaces} As a byproduct of the microlocal techniques used in the proof of Theorem \ref{mainthm1}, we obtain an improvement over a result of Flaminio \cite{Flaminio-92} on flow-invariant distributions on convex surfaces. We now assume that $g$ is a positively-curved metric on $M := \Ss^2$.
It was proved by Flaminio \cite{Flaminio-92} that any even flow-invariant distribution $u \in \mc{D}'(SM)$ (that is, whose expansion $u = \sum_{k \in \Z} u_k$ only contains even Fourier modes) on the unit tangent bundle $SM$ is determined by its push forward $\pi_* u \in \mc{D}'(M)$ (where $\pi : SM \to M$ denotes the projection). This statement is no longer true in higher dimensions as proved by Sutton in \cite{Sutton-03}. We prove the following stability estimate, complementing \cite{Flaminio-92}:

\begin{theorem}
\label{theorem:sphere-intro}
Assume that $(M,g)$ has positive curvature. Then for all $s \in \R$, there exists a constant $C > 0$ such that the following inequality holds: for all $u \in C^\infty(SM)$,
\begin{equation}
\label{equation:stabilite1}
\|u\|_{H^s} \leq C \left( \|Xu\|_{H^s} + \|u_0\|_{H^{s+1}}+\|u_1\|_{H^{s+1}}\right).
\end{equation}
Moreover, if $u \in \mc{D}'(SM)$ is a distribution, then for all $s \in \R$,
\[
Xu \in H^s(SM), u_0,u_1 \in H^{s+1}(SM) \implies u \in H^s(SM),
\]
and \eqref{equation:stabilite1} holds. In particular, if $u \in \mc{D}'(SM)$, $Xu=0$, $u_0+u_1 = 0$, then $u=0$.
\end{theorem}

We actually prove a more general result for the operator $X + \pi^*\phi$ (and not just $X$), see Theorem \ref{theorem:sphere} below. The assumption could be relaxed by simply requiring the existence of conjugate points along \emph{every geodesic}. Theorem \ref{theorem:sphere-intro} can be seen as a quantitative improvement over \cite{Flaminio-92}. We also formulate some remarks relative to Theorem \ref{theorem:sphere-intro}:

\begin{enumerate}[label=(\roman*)]
\item An immediate consequence of Theorem \ref{theorem:sphere-intro} is the following surprising result: if $u \in \mc{D}'(SM)$, $Xu=0$ and $u_0+u_1 \in C^\infty(M,\Omega_0\oplus\Omega_1)$, then $u \in C^\infty(SM)$. In other words, the regularity of the first two Fourier modes of a flow-invariant distribution also determines the regularity of the distribution.

\item If $u \in \mc{D}'(SM)$ is flow-invariant, namely, $X u = 0$, and $u \in H^s(SM)$, then the averaging result of Gérard-Golse \cite[Theorem 2.1]{Gerard-Golse-92} for kinetic equations implies that $u_0, u_1 \in H^{s+1/2}(SM)$. Theorem \ref{theorem:sphere-intro} can thus be seen as a converse to this statement: if $u \in \mc{D}'(SM)$ is flow-invariant and $u_0, u_1 \in H^{s+1}(SM)$, then $u \in H^s(SM)$.
\end{enumerate}

\subsection{Holomorphic line bundles over twistor space} One of the main results of \cite{Bohr-Paternain-21} was the classification of holomorphic vector bundles over the twistor space of simple Riemannian surfaces: there, a `transport Oka--Grauert principle' was proved, showing that all holomorphic vector bundles are actually trivial. In our last result, we classify all holomorphic line bundles over $Z$ having certain distributional behaviour on the boundary $\partial Z$, and up to isomorphism in $Z^\circ$, see Theorem \ref{modulianosov}. To keep this introductory discussion accessible, we restrict here to topologically trivial holomorphic line bundles extending smoothly up to the boundary $\partial Z$ and refer the interested reader to \S\ref{section:line-bundles} for more general results.

If $U \subset Z$ is an open subset, we denote by $\Omega^{p,q}(U)$ the set of smooth $(p,q)$-forms on $U$. This is well defined also if $U\cap \partial Z \neq \emptyset$, though in this case  $\Omega^{p,q}(U) \not \subset \Omega^{p+q}(U)$ due to the degeneracy of the complex structure -- see \textsection \ref{section:dolbeaut} for details.
We think of $\Omega^{p,q}(Z)$ as the space of $(p,q)$-forms that are smooth up to the boundary $\partial Z$. We let $\eps := Z \times \C \to Z$ be the trivial line bundle over $Z$. A holomorphic structure on $\eps$, smooth up to the boundary, is the data of a \emph{partial connection}
\[
\bar \partial + \alpha : \Omega^0(Z) \to \Omega^{0,1}(Z)
\]
where $\alpha \in \Omega^{0,1}(Z)$ satisfies the integrability condition $\bar \partial \alpha = 0$. Two holomorphic structures are equivalent on a subset $U\subset Z$, written 
$(\eps, \bar \partial + \alpha)\simeq (\eps,\bar \partial+\tilde{\alpha})$ on $U$, if there exists a function $\varphi \in C^\infty(U,\C^\times)$ with $\alpha = \varphi^{-1} (\bar \partial + \tilde{\alpha})\varphi$.

There is a natural fibration $\fib\colon SM \times \D \to Z$ (where $\D$ denotes the closed unit disk) whose fibres correspond to the $S^1$-orbits of a smooth vector field $\mathbf{V}$, see \S\ref{section:twistor}. In order to study objects on $Z$, it is often more convenient to lift them on $SM \times \D$ as $\mathbf{V}$-invariant objects. We will prove the following classification result:

\begin{theorem}\label{modulianosov-intro} Suppose $Z$ is the twistor space of an Anosov surface.

\begin{enumerate}[label=\emph{(\roman*)}]
\item There exists $\tau \in  \Omega^{1}(SM \times \D)$ such that the following holds. Given $a\in C^\infty(SM)$ with $a_k =0$ for $k<-1$, we form
\begin{equation}
\label{equation:alpha-intro}
\alpha = \fib_*\left(\sum_{k\ge -1} a_k \omega^{k+1} \tau\right)\in \Omega^{0,1}(Z).
\end{equation}
Then the partial connection $\bar \partial + \alpha$ is a holomorphic structure on $\eps$, that is, $\bar \partial \alpha = 0$. Conversely, any holomorphic structure on $\eps$ is equivalent to $\bar \partial + \alpha$ for some $\alpha \in \Omega^{0,1}(Z)$ of the form \eqref{equation:alpha-intro}.
\smallskip
\item  Let $\iota_M\colon M\rightarrow Z$ be the inclusion as $0$-section and consider $\alpha \in \Omega^{0,1}(Z)$ as in \eqref{equation:alpha-intro} and
\[
 \beta = \fib_*\left(\sum_{k\ge 0} a_{2k} \omega^{2k} \bar \tau \wedge \tau\right) \in \Omega^{1,1}(Z).
 \]
Then
\begin{equation}
\label{equation:pouet}
\def\arraystretch{1.5}
\left \{ \begin{array}{ll}
(i) & [\iota_M^*\alpha]=0 \in  H^{0,1}_{\bar \partial}(M)/H^1(M,\Z)\\
(ii) & [\iota_M^*\beta] =0 \in H^2(M)
\end{array} \right. \Longleftrightarrow (\eps, \bar \partial + \alpha) \simeq (\eps, \bar \partial) \text{ on } Z^\circ.
\end{equation}
\end{enumerate}
\end{theorem}

\medskip

Here $H^{0,1}_{\bar \partial}(M)/H^1(M,\Z) = \mathrm{Pic}_0(M) \simeq \C^g/\Z^{2g}$
 is the identity component of the Picard group $\mathrm{Pic}(M)$ of the Riemann surface $(M,[g])$, that is, of the group of holomorphic line bundles over $M$, see \S\ref{section:linebundlesonM} where this is recalled. The theorem is in line with the transport Oka--Grauert principle for simple surfaces \cite{Bohr-Paternain-21}, where $H^{0,1}_{\bar \partial}(M)=H^2(M)=0$, such that the two obstructions disappear.

 The equivalence \eqref{equation:pouet} shows that $\mathrm{Pic}_0(Z,Z^\circ)$, the space of holomorphic structures over the trivial line bundle that are smooth up to the boundary and considered up to isomorphism in the interior of the twistor space $Z^\circ$, is given by 
\[
\mathrm{Pic}_0(Z,Z^\circ) \simeq H^{0,1}_{\bar \partial}(M)/H^1(M,\Z) \times H^2(M) = \C^g/\Z^{2g} \times \C.
\]
The results of Theorem \ref{modulianosov-intro} are generalized in Theorem \ref{modulianosov} below in order to capture more singular holomorphic structures which may blow up (polynomially) at $\partial Z$.

\subsection{Organisation of the article} In \S\ref{section:unit}, we recall standard results on the geometry, dynamics and Fourier analysis of the unit tangent bundle of a Riemannian surface. In \S\ref{section:twistor}, we recall the construction from \cite{Bohr-Paternain-21} of the transport twistor space $Z$. Section \S\ref{section:invariant} contains the proof of the twistor correspondence for invariant distributions, namely, Theorem \ref{mainthm1}. In \S\ref{section:intermezzo}, as a byproduct of the microlocal techniques of \S\ref{section:invariant}, we prove a result generalizing earlier work of Flaminio \cite{Flaminio-92} on the regularity of invariant distributions on convex surfaces (see Theorem \ref{theorem:sphere}). Finally, we classify holomorphic line bundles over twistor space in \S\ref{section:line-bundles}. Some standard results of microlocal analysis are recalled in Appendix \ref{appendix} for the reader's convenience (elliptic operators, propagation of singularities). 

\subsection*{Acknowledgements}  We are grateful to Colin Guillarmou for very helpful discussions related to Theorem \ref{mainthm1}.
We are also grateful to the referees for their suggestions for corrections and improvements.

\section{Geometry and analysis on the unit tangent bundle}

\label{section:unit}

\subsection{Unit tangent bundle}

In this paragraph, we review some basics of Riemannian geometry on surfaces. We refer to \cite{Paternain-Salo-Uhlmann-book} for a more extensive treatment.

\subsubsection{Geometry of the unit tangent bundle} \label{section:geometrySM}

Let $(M,g)$ be a closed oriented Riemannian surface. The unit tangent bundle $\pi : SM \to M$ is defined as
\[
SM := \left\{ (x,v) \in TM ~|~ |v|_g = 1\right\}.
\]
Let $X \in C^\infty(SM,T(SM))$ be the generator of the geodesic flow $(\varphi_t)_{t \in \R}$, and let $V \in C^\infty(SM,T(SM))$ be the generator of the $\mathrm{SO}(2)$-action in the circle fibres. Define $H := [V,X]$. The frame $\left\{X,H,V\right\}$ is an orthonormal frame with respect to the Sasaki metric on $SM$. We introduce the vertical distribution $\V := \R V$, and the horizontal one $\HH := \R H$. (Note that, in the usual terminology, the horizontal distribution also includes the span of $X$.) The following commutation formulas hold:
\begin{equation}
\label{equation:commutators}
[X,H] = K V, \qquad [V,X]= H, \qquad [H,V] = X,
\end{equation}
where $K$ is the sectional curvature. The raising and lowering operators $\eta_\pm$ are defined as
\begin{equation}
\label{equation:etapm}
\eta_\pm := \dfrac{1}{2}\left(X \mp iH \right).
\end{equation}
It is immediate to check from \eqref{equation:etapm} that $\eta_+^*=-\eta_-$, where the formal adjoint is computed with respect to the $L^2$ scalar product on $SM$ induced by the Liouville measure.
Combining \eqref{equation:commutators} and \eqref{equation:etapm}, we also obtain the following commutation formulas:
\begin{equation}
\label{equation:commutateur2}
[-iV,\eta_+] = \eta_+, \qquad [-iV,\eta_-]=-\eta_-, \qquad [\eta_+,\eta_-] = \dfrac{i K V}{2}.
\end{equation}

We now introduce the dual distributions on $T^*(SM)$. Let $\alpha$ be the Liouville contact $1$-form defined by $\iota_X \alpha = 1, \iota_X d\alpha = 0$. It satisfies $\alpha(\HH\oplus\V)=0$. We introduce $\beta,\gamma \in C^\infty(SM,T^*(SM))$ such that 
\[
\beta(H)=1, \beta(X)=\beta(V)=0, \qquad \gamma(V)=1,\gamma(X)=\gamma(H)=0.
\]
We set $E_0^* := \R \alpha$. The cohorizontal distribution is defined as $\HH^*:=\R\beta$, whereas the covertical is $\V^*:=\R\gamma$. Finally, we define the characteristic set $\Sigma := \HH^* \oplus \V^*$.

\subsubsection{Fourier decomposition}\label{section:fourier}

Any function $f \in C^\infty(SM)$ can be decomposed in Fourier modes by simply freezing the basepoint $x \in M$ and considering the Fourier decomposition on $S_xM \simeq \Ss^1$. More precisely, any $u\in C^\infty(SM)$ can be expanded as:
\begin{equation}
\label{equation:decomposition}
u = \sum_{k \in \Z}u_{k}\;\;\text{with}\;u_{k}\in C^\infty(M,\Omega_k), \qquad \Omega_k(x) := \ker(V-ik)|_{C^\infty(S_xM)}.
\end{equation}
Here, $\Omega_k \to M$ is a holomorphic line bundle over $M$ and there is an implicit identification
\begin{equation}
\label{equation:identification}
\pi_k^* : C^\infty(M,\Omega_k) \to C^\infty(SM)
\end{equation}
of a section of $\Omega_k$ with a function $SM$ annihilated  by $V-ik$ . The dual operator to $\pi_k^*$ is denoted by ${\pi_k}_*$ and the $L^2$ orthogonal projection onto the $k$-th Fourier mode is then simply given by $\pi_k^* {\pi_k}_*$. 

The decomposition \eqref{equation:decomposition} is orthogonal with respect to the $L^2$ scalar product on $SM$. Moreover, $\Omega_k$ can be naturally identified with $\kappa^{\otimes k}$, the $k$-th power of the canonical line bundle. 
More generally, the expansion \eqref{equation:decomposition} holds for any distribution $u \in \mc{D}'(SM)$ by writing $u = \sum_k u_k$, where $u_k \in \mc{D}'(M,\Omega_k) \subset \mc{D}'(SM)$ is defined by setting:
\[
\forall \varphi \in C^\infty(SM), \qquad \langle u_k,\varphi \rangle := \langle u, \varphi_k \rangle.
\]
A function/distribution is said to be even (resp.\,odd) if its Fourier decomposition only contains even (resp.\,odd) Fourier modes.

The following identity is known as the \emph{localized Pestov identity}:

\begin{lemma}
For all $u \in C^\infty(M,\Omega_k)$:
\begin{equation}
\label{equation:pestov}
\|\eta_+u\|^2_{L^2(SM)} =  \|\eta_-u\|^2_{L^2(SM)} - \frac{k}{2} \langle Ku,u\rangle_{L^2(SM)}.
\end{equation}
\end{lemma}

The proof is a straightforward computation using \eqref{equation:commutateur2}.
The operators $\eta_\pm$ act as lowering/raising operators on this decomposition, namely
\[
\eta_\pm : C^\infty(M,\Omega_k) \to C^\infty(M,\Omega_{k \pm 1}).
\]
Moreover, $\eta_\pm$ are elliptic differential operators of order $1$ (when restricted to sections of $\Omega_k$). The operator $\eta_-$ (resp. $\eta_+$) acting on $\Omega_k$ for $k \geq 0$ is conjugate to the $\overline{\partial}$ operator (resp. $\partial$ operator) acting on $\kappa^{\otimes k}$, see \cite[Section 2]{Paternain-Salo-Uhlmann-14-2}.  In particular, the space of holomorphic $m$-differentials on $M$ is given by
\begin{equation}\label{def:mdifferential}
\mathcal H_m = C^\infty(M,\kappa^{\otimes m})\cap \ker \bar \partial \simeq C^\infty(M,\Omega_m)\cap \ker \eta_-.
\end{equation}
Since $\partial$ and $\bar \partial$ only depend on the conformal structure of $(M,g)$, the dimensions of the kernels of $\eta_\pm$ are conformally invariant. They can be easily computed using the Pestov identity \eqref{equation:pestov} and the Riemann-Roch theorem.

\begin{lemma}
\label{lemma:injectivity-eta}
Let $M$ be a closed oriented surface of genus $g \geq 0$. Then, the following holds:
\begin{enumerate}[label=\emph{(\roman*)}]
\item If $g=0$, then for all $k \in \Z$,
	\begin{itemize}
		\item $\dim \ker \eta_+|_{C^\infty(M,\Omega_k)} = (2k+1)_+$,
		\item $\dim \ker \eta_-|_{C^\infty(M,\Omega_k)}=(1-2k)_+$.
	\end{itemize}
\item If $g=1$, then for all $k \in \Z$, $\ker \eta_+|_{C^\infty(M,\Omega_k)}= \ker \eta_-|_{C^\infty(M,\Omega_k)}=\C$.
\item If $g \geq 2$, then:
	\begin{itemize}
		\item Kernels of $\eta_+$: For $k \leq-2$, $\dim \ker \eta_+ = -(2k+1)(g-1)$; for $k=-1$, $\dim \ker \eta_+ = g$; for $k=0$, $\dim \ker \eta_+=1$; for $k \geq 1$, $\eta_+$ is injective.
\item Kernels of $\eta_-$: For $k \geq 2$, $\dim \ker \eta_- = (2k-1)(g-1)$; for $k=1$, $\dim \ker \eta_- = g$; for $k =0$, $\dim \ker \eta_- = 1$; for $k \leq -1$, $\eta_-$ is injective.
	\end{itemize}
\end{enumerate}
\end{lemma}

We refer to \cite[Lemma 2.1]{Paternain-Salo-Uhlmann-14-2} for a proof.
We now introduce the following terminology:

\begin{definition}\label{def:fibrehol}
A function $u \in C^\infty(SM)$ is said to be \emph{fibrewise holomorphic} (resp. \emph{antiholomorphic}) if $u_k = 0$ for all $k < 0$ (resp. $u_k = 0$ for all $k > 0$).
\end{definition}

The Szeg\H{o} projector $\mc{S} : L^2(SM) \to L^2(SM)$ is the $L^2$-orthogonal projector onto holomorphic functions, namely $\mc{S}(\sum_k u_k) = \sum_{k \geq 0} u_k$.
We have the commutation formula:
\begin{equation}
\label{equation:szego-commutation}
[\mc{S},X]f =\eta_+ f_{-1} - \eta_- f_0.
\end{equation}

\subsection{Projective dynamics on the characteristic set}

The geodesic flow $(\varphi_t)_{t \in \R}$ on $SM$ lifts naturally to $T(SM)$ and $T^*(SM)$ by defining for $v \in SM$, $w \in T_v(SM)$ and $\xi \in T^*_v(SM)$:
\[
\Phi^{T(SM)}_t(v,w) := (\varphi_t(v), d\varphi_t(v)w), \quad \Phi^{T^*(SM)}_t(v,\xi) := (\varphi_t(v), d\varphi_t^{-\top}(v)\xi).
\]
In the following, we will be mostly interested in the flow on the cotangent bundle. By the definition of $d\varphi_t^{-\top}(v)$ we have
\[\eta(t):=d\varphi_t^{-\top}(v)\xi\in T^{*}_{\varphi_{t}(v)}(SM)\]
if and only if
\begin{equation}
\xi=\eta(t)\circ d\varphi_t(v).
\label{eq:eta}
\end{equation}
The Hamiltonian flow of $\langle \xi,X(v)\rangle$ is precisely $(\Phi^{T^*(SM)}_t)_{t \in \R}$  and hence it
preserves $\Sigma$.
In what follows it is convenient to consider the projectivised action of $(\Phi^{T^*(SM)}_t)_{t \in \R}$ on $\Sigma$. Note that $\Sigma\subset T^*(SM)$ is a sub-bundle with fibres
\begin{equation*}
\Sigma(v) =\{\xi\in T^*_v(SM): \langle \xi,X(v)\rangle =0\},\quad v\in SM.
\end{equation*}
We let $\Lambda(SM)$ be the 4-manifold obtained by projectivising $\Sigma(v)$ for all $v\in SM$. In other words $\Lambda(SM)$ is the bundle over $SM$ whose fibre over $v$ consists of all 1-dimensional subspaces of $\Sigma(v)$. Note that $\V^*$ and $\HH^*$ are sections of this bundle.

Finally observe that the flow $(\Phi^{T^*(SM)}_t)_{t \in \R}$ naturally induces a flow on $\Lambda(SM)$ that we denote by $(\Phi^{\Lambda(SM)}_t)_{t \in \R}$.

\subsubsection{Twist property of the vertical bundle}

\label{sssection:twist}

Consider the hypersurface $\Lambda_{\HH^*}\subset \Lambda(SM)$ given by $\HH^*(SM)$.
The following property of the geodesic flow is well-known:



\begin{proposition} The flow $(\Phi^{\Lambda(SM)}_t)_{t \in \R}$ is transversal to $\Lambda_{\HH^*}$.
\label{prop:coV}
\end{proposition}

\begin{proof} A proof of this proposition is given in \cite[Proposition 2.40]{Paternain-99} for the induced flow acting on the tangent bundle of $SM$.
We provide a direct self-contained proof for the dual statement on the cotangent bundle of $SM$.
Take any $\xi\in \Sigma(v)$ and let $\eta(t)=d\varphi_t^{-\top}(v)\xi$.
Write $\eta(t)=x(t)\gamma+y(t)\beta$ for some functions $x(t)$ and $y(t)$ and use \eqref{eq:eta} to derive
\[\xi=x(t)\gamma\circ d\varphi_{t}(v)+y(t)\beta\circ d\varphi_{t}(v)=x(t)\varphi_{t}^*\gamma+y(t)\varphi_{t}^*\beta.\]
Differentiating this identity with respect to $t$ and using the definition of the Lie derivative ${\mathcal L}_{X}$ we obtain
\[0=\dot{x}\,\varphi_{t}^*\gamma+x \,\varphi_{t}^*({\mathcal L}_{X}\gamma)+\dot{y}\,\varphi_{t}^*\beta+y\,\varphi^*_{t}({\mathcal L}_{X}\beta).\]
Using the structure equations \eqref{equation:commutators} we see that ${\mathcal L}_{X}\beta=\gamma$ and ${\mathcal L}_{X}\gamma=-K\beta$.
Hence
\[0=(\dot{y}-Kx)\varphi_{t}^*\beta+(y+\dot{x})\varphi_{t}^*\gamma.\]
Since $\{\beta,\gamma\}$ are linearly independent we have the pair of equations 
\begin{equation}
\dot{y}-Kx=0,\;\;\;\;y+\dot{x}=0.
\label{eq:jacobico}
\end{equation}
The second equation gives the transversality property right away as follows. Consider the function
$h:\Sigma\to\R$ given by $h(v,\xi)=\xi(V(v))$. Clearly 
$$h^{-1}(0)=\{(v,\,\HH^*(v)):\;v\in SM\}$$
and
\[h\left (\Phi_{t}^{T^*(SM)}(v,\xi)\right)=\eta(t)(V(\varphi_{t}(v)))=x(t).\]
Taking the derivative at $t=0$ gives
\[dh_{(v,\xi)}(\mathbf{X})=\dot{x}(0)=-y(0)=-\xi(H(v)),\]
where $\mathbf{X}$ is the infinitesimal generator of $(\Phi^{T^*(SM)}_t)_{t \in \R}$. Hence, we deduce that if $(v,\xi)\in h^{-1}(0)$ and
$\xi\neq 0$, then $dh_{(v,\xi)}(\mathbf{X})\neq 0$ and the transversality claim follows.


\end{proof}

\begin{remark} The same transversality property applies to $\V^*$ provided the Gaussian curvature has a sign.
Note also that the proof indicates that $\HH^*$ moves in the counterclockwise direction as we apply the flow as depicted in Figure \ref{figure:cones}.

\end{remark}

\subsubsection{Positive curvature}

When the Gaussian curvature $K>0$, every geodesic will have conjugate points. In this instance the dynamics of the projectivised flow $(\Phi^{\Lambda(SM)}_t)_{t \in \R}$ is easy to visualise: when we flow $\HH^*$ (forward and backwards in time) it will turn around the circle
indefinitely and the flow will meet the two hypersurfaces $\Lambda_{\HH^*}$ and $\Lambda_{\V^*}$ transversally infinitely many times.

\label{sssection:projective-positive}

\subsubsection{Non conjugate points and nonpositive curvature}

\label{sssection:projective-negative}

\begin{center}
\begin{figure}[htbp!]
\includegraphics{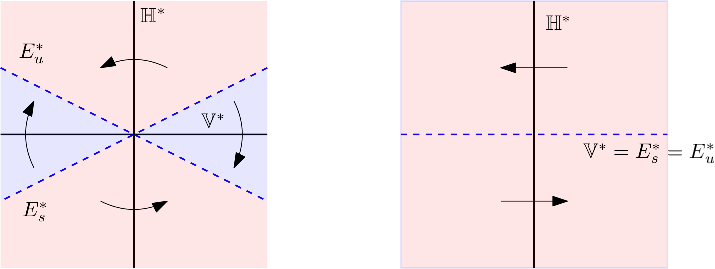}
\caption{The lifted dynamics on $\Sigma$. Left picture: dynamics near a point where $E_s^*$ and $E_u^*$ are transverse; the light blue region represents the cone $\mc{C}$. Right picture: dynamics over the set where $E_s^*=E_u^*=\V^*$; the cone $\mc{C}$ collapses here to a line.}
\label{figure:cones}
\end{figure}
\end{center}

The proof of Proposition \label{prop:coV} gives the Jacobi equation $\ddot{y}+Ky=0$ for the component $y$ of $\eta(t)$. The surface is said to have no conjugate points if, for any geodesic, any non-trivial solution $y$ of the Jacobi equation with $y(0) = 0$ vanishes only at $t = 0$. Equivalently, $d\varphi_{t}^{-\top}(\HH^*(v))\cap\HH^*(\varphi_t(v))=\{0\}$ for all $t\neq 0$ and all $v\in SM$. Under the no conjugate points assumption, \cite{Hopf_48} gives the existence of two limit solutions of the Riccati equation $\dot{r}+r^2+K=0$ along a geodesic that are defined for all $t\in\R$. The existence of these two solutions translate into the existence of two invariant subbundles that we denote by $E^*_{u}$ and $E^{*}_{s}$ obtained by the following limiting procedure:
\[E^*_{u}(v)=\lim_{t\to\infty}d\varphi_{t}^{-\top}(\HH^*(\varphi_{-t}(v))\]
\[E^*_{s}(v)=\lim_{t\to\infty}d\varphi_{-t}^{-\top}(\HH^*(\varphi_{t}(v))\]
When $K\leq 0$, the bundles are continuous in $v$. Moreover, a classical result of Eberlein \cite{Eberlein_73} guarantees that the geodesic flow is Anosov if and only if $E^{*}_{u}(v)\cap E^{*}_{s}(v)=\{0\}$ for all $v\in SM$.  Since the bundles $E^{*}_{u}$ and $E^{*}_{u}$ are invariant they determine uniquely an invariant cone $\mc{C} \subset T^*(SM)$ defined as the cone delimited by $E^{*}_{u}$ and $E^{*}_{u}$ but not containing $\HH^*$, see Figure \ref{figure:cones}. The cone could collapse to a line as in the case of a geodesic completely contained in a region of zero curvature.


\section{Twistor space of surfaces}

\label{section:twistor}

\subsection{Definition}

In this section, we review the twistor space construction from \cite{Bohr-Paternain-21}. Here $(M,g)$ may be any oriented Riemannian surface, possibly with non-empty boundary $\partial M$. As a smooth manifold, the twistor space equals the total space of the unit disk bundle,
\begin{equation}
Z =\{(x,v)\in TM: g(v,v)\le 1\}.
\end{equation}
The complex structure is encoded in a complex $2$-plane distribution $\mathscr D\subset T_\C Z=TZ\otimes \C$ which is {\it involutive} (that is, $[\mathscr D,\mathscr D]\subset\mathscr D$) and satisfies the {\it (non-)degeneracy condition}
\begin{equation}\label{Dproperties}
\mathscr D\cap \overline{\mathscr D} = 0 \text{ on } Z\backslash SM\quad\text{ and } \quad \mathscr D\cap \overline{\mathscr D} = \spn_\C (X)\text{ on } SM.
\end{equation}
In order to define $\mathscr D$, it is convenient to introduce the fibration 
\[ \fib \colon SM\times \D\rightarrow Z,\;\;\; \fib (x,v,\omega)=(x,\omega v),\]
where the product $\omega v\in T_x M$ is explained by the complex structure that $g$ and the orientation induce on $M$. In fact, the fibres of $\fib$ are given by the integral curves of the following vector field on $SM\times \D$:
\begin{equation}
\mathbf{V} = V +  i (\bar \omega \partial_{\bar \omega}  - \omega \partial_\omega)
\end{equation}
Using the commutator relations \eqref{equation:commutateur2} one shows that the complex vector fields $\omega^2\eta_++\eta_-$ and $\partial_{\bar \omega}$ lie in the eigenspace $\ker (\mathcal L_{\mathbf{V}} + i)$, where $\mathcal L_{\mathbf{V}}$ denotes the Lie derivative. Hence the complex line bundles spanned by these vector fields are $S^1$-invariant and the push-forward
\begin{equation}
\mathscr D =\fib_*\spn_\C \left(\omega^2\eta_++\eta_-,\partial_{\bar \omega}\right)\subset T_\C Z
\end{equation}
is a well defined $2$-plane distribution on $Z$. For more details, including the proof that $\mathscr D$ is involutive and satisfies \eqref{Dproperties}, we refer to \cite[Lemma 4.1]{Bohr-Paternain-21}. The relation of $Z$ with the classical twistor spaces for projective structures is discussed in  \cite[\textsection 4.4]{Bohr-Paternain-21}

In the interior $Z^\circ$ of the twistor space, the non-degeneracy condition \eqref{Dproperties} implies that $\mathscr D$  is the $(-i)$-eigenspace of an integrable complex structure. Hence $Z^\circ$ is a classical complex surface with $(0,1)$-tangent bundle given by
\begin{equation}
T^{0,1}Z^\circ  = \mathscr D\vert_{Z^\circ}.
\end{equation}
If $U\subset Z$ is an open set with $U\cap \partial Z = \emptyset$, then a smooth function $f\colon U\rightarrow \C$ is holomorphic if and only if $df\vert_{\mathscr D} =0$ on $U$. Equivalently, the pull-back $h=\fib^*f\in C^\infty(\fib^{-1}(U))$ satisfies the Cauchy--Riemann equations
\begin{equation}
(\omega^2\eta_+  + \eta_-)h = 0 \quad \text { and } \quad \partial_{\bar \omega} h=0.
\end{equation}
Note that the latter two characterisations of holomorphicity also make sense in case $U\cap \partial Z\neq \emptyset$; in particular the definition $\mathcal A(U)=\{f\in C^\infty(U): df \vert_{\mathscr D}=0\}$ encapsulates the definitions for $\mathcal A(Z)$ and $\mathcal A(Z^\circ)$ given in the introduction. Similarly, other definitions from complex geometry carry over to twistor space $Z$ in a way that is {\it smooth up to the boundary}---we will discuss further instances of this below.

\subsection{Dolbeaut complex}\label{section:dolbeaut} Let us revise the Dolbeaut complex
\begin{equation}\label{dolbeautcx}
\Omega^{p,0}(U)\xrightarrow{\bar \partial} \Omega^{p,1}(U)\xrightarrow{\bar \partial} \Omega^{p,2}(U),\qquad p\in\{0,1,2\}.
\end{equation}
On an open set $U\subset Z$ with $U\cap \partial Z= \emptyset$, where the complex structure does not degenerate, this is defined as usual in complex geometry.
We now give a description  in terms of the fibration $\fib\colon SM \times  \D\rightarrow Z$, which allows to extend to the case $U\cap \partial Z\neq \emptyset$. For $p=0$ this was already done in \cite{Bohr-Paternain-21}; however, here we are interested also in the case $p=1$.

On $SM\times \D^\circ$ we define complex $1$-forms $\tau$ and $\gamma$ by the relations
\begin{equation}\label{taugamma}
\begin{array}{l}
\tau(\omega^2\eta_++\eta_-)=\gamma(\partial_{\bar \omega})=1\\
 \tau(\partial_{\bar \omega})=\gamma(\omega^2\eta_++\eta_-)=0
\end{array}
,\qquad \tau,\gamma = 0 \text{ on } \spn_\C(\bar \omega^2\eta_-+\eta_+,\partial_\omega,\mathbf V).
\end{equation}
While this forces $\tau$ itself to blow up as $\vert \omega \vert \rightarrow 1$, both $(1-\vert \omega\vert^4)\tau$ and $\gamma$ extend smoothly to $SM\times \D$, see also \cite[Lemma 4.2]{Bohr-Paternain-21}. There it was also shown that
\begin{equation}\label{lveigen1}
\tau,\gamma \in \ker(\mathcal L_{\mathbf{V}} -i),\quad \bar \tau,\bar \gamma\in \ker (\mathcal L_{\mathbf{V}}+i).
\end{equation}
In order to express differential forms on $Z$ in this co-frame we consider, for $U\subset Z$ open, the coefficient spaces
\begin{equation*}
\mathcal I^m(U)=\{h\in C^\infty(\fib^{-1}(U)):(\mathbb V - i m)h=0\},\quad  m\in \Z.
\end{equation*}
Note that if $U\cap \partial Z\neq \emptyset$, this space contains functions that are smooth up to the boundary of $Z$.
Away from the boundary, if $U\cap \partial Z=\emptyset$, we define maps for all $p,q\in \{0,1,2\}$ as follows:
\begin{equation}\label{phipq}
\Phi^{p,q} \colon \underbrace{\mathcal I^{p-q}(U)\times \dots \times \mathcal I^{p-q}(U)}_{N\text{-times}, ~N= {2 \choose p}\cdot {2 \choose q}} \rightarrow \Omega^{p,q}(U),\quad h=(h_1,\dots,h_N)\mapsto \fib_*(\varepsilon_h^{p,q}),
\end{equation}
with $\varepsilon_h^{p,q}\in \Omega^{p+q}(\fib^{-1}(U))$  as detailed in Table \ref{table}. For example, if $p=q=1$, then $N=4$ and $h=(h_1,h_2,h_3,h_4)\in \mathcal I^0(U)^4$ is sent to
\begin{equation}
\Phi^{1,1}(h) = \fib_*(h_1 \bar \tau\wedge \tau + h_2 \bar \tau\wedge \gamma + h_3 \bar \gamma \wedge \tau + h_4 \bar \gamma \wedge  \gamma).
\end{equation}
Here the push-forward is well defined because of \eqref{lveigen1} and the observation $\ker(\mathcal L_\mathbf V -i m_1)\wedge \ker(\mathcal L_\mathbf V -i m_2)\subset \ker(\mathcal L_\mathbf V -i (m_1+m_2))$, which follows from the Leibniz rule. We further define operators
\begin{equation}\label{dpq}
D^{p,q}\colon \mathcal I^{p-q}(U)^{{2\choose p}{2\choose q}} \rightarrow \mathcal I^{p-q-1}(U)^{{2\choose p}{2\choose {q+1}}},\quad h\mapsto D^{p,q}h =\delta^{p,q}_h,
\end{equation}
with $\delta^{p,q}_h$ also given in Table \ref{table}.

\begin{table}
\begin{tabular}{c|c|c|c}
$\varepsilon_h^{p,q}$ & $q=0$ & $q=1$ & $q=2$\\[.5ex]
 \hline 
 \,&&&\\[-2ex]
 $p=0$ & $h_1$ & $h_1\tau + h_2 \gamma$&  $h_1 \tau\wedge \gamma$  \\[1ex]
 $p=1$ & $h_1 \bar \tau + h_2 \bar \gamma$ & \scalebox{.8}{$\left(\begin{array}{c}
 \quad h_1 \bar \tau\wedge \tau + h_2 \bar \tau\wedge \gamma \\
 + h_3 \bar \gamma \wedge \tau + h_4 \bar \gamma \wedge  \gamma
 \end{array}\right)$}
 & 
 \scalebox{.8}{$\left(\begin{array}{l}
 \quad h_1 \bar \tau \wedge \tau \wedge \gamma\\
 \,+ h_2 \bar \gamma \wedge \tau \wedge \gamma
 \end{array}\right)$}
  \\[2ex]
 $p=2$ & $h_1 \bar \tau \wedge \bar \gamma$ & $h_1 \bar \tau \wedge \bar \gamma \wedge \tau +h_2 \bar \tau \wedge \bar \gamma \wedge \gamma $&  $h_1 \bar \tau\wedge \bar \gamma \wedge \tau\wedge \gamma$
\end{tabular}
\vspace{2em}
\begin{tabular}{l}
\,
\end{tabular}\\
\begin{tabular}{c|c|c}
$\delta_h^{p,q}$ & $q=0$ & $q=1$
\\[.5ex]
 \hline 
  \,&&\\[-1.75ex]
 $p=0$ & $\big(\xi h_1,\partial_{\bar \omega}h_1\big)$ & $-\partial_{\bar \omega} h_1 + \xi h_2$\\[1.5ex]
 $p=1$ & 
\scalebox{.8}{$\left(\begin{array}{ll}
		-\xi h_1 + \frac{\omega}{2}(1-|\omega|^4)Kh_2, & -\partial_{\bar \omega} h_1 - \frac{2\omega |\omega|^2}{1-|\omega|^4}h_1,\\[1.5ex]
		-\xi h_2 - \frac{2\omega}{1-|\omega|^4} h_1, & -\partial_{\bar \omega} h_2 
	\end{array}\right)
	
 $} & 
 \scalebox{.8}{$\left( \begin{array}{l}
\left( \partial_{\bar \omega}+ \frac{2\omega|\omega|^2}{1-|\omega|^4}\right) h_1-\xi h_2 + \frac  {\omega(1-|\omega|^4)}2 K h_4,\\[1.5ex]
 \partial_{\bar \omega}h_3-\xi h_4 - \frac{2\omega}{1-|\omega|^4} h_2
 \end{array}
 \right)$
 }
  \\[2.5ex]
 $p=2$ & $\left(\xi h_1,\big(\partial_{\bar \omega}+\frac{2\omega|\omega|^2}{1-|\omega|^4}\big)h_1\right)$ & $-\partial_{\bar \omega}h_1+ \xi h_2 - \frac{2\omega|\omega|^2}{1-|\omega|^4}h_1$
\end{tabular}
\begin{tabular}{l}
\,
\end{tabular}\\
\caption{Definition of $\varepsilon_h^{p,q}$ from \eqref{phipq} and $\delta^{p,q}_h$ from \eqref{dpq}. Here $\xi = \omega^2\eta_++\eta_-.$}
\label{table}
\end{table}

\begin{proposition}\label{intertwine}
Let $U\subset Z$ be an open set with $U\cap \partial Z = \emptyset$. Then the map $\Phi^{p,q}$ defined in \eqref{phipq} is an isomorphism and intertwines the Dolbeaut complex from \eqref{dolbeautcx} with the operators from \eqref{dpq}, that is, for all $0\le p \le 2$ and $0\le q \le 1$
 we have
 \begin{equation}
 \bar \partial \circ \Phi^{p,q} = \Phi^{p,q+1}\circ D^{p,q}.
 \end{equation}
 
 \end{proposition}

\begin{proof}
Let $k=p+q$, then a $k$-form $\alpha \in \Omega^k(U)$ has bi-degree $(p,q)$ if and only if 
$\alpha(\xi_1,\dots,\xi_{k})=0$ for any collection of vector fields with $\xi_1,\dots,\xi_{p'}\in T^{1,0}U$ and $\xi_{p'},\dots,\xi_k\in T^{0,1}U\equiv \mathscr D\vert_U$ for some $p'\neq p$. This is further equivalent to the condition
\begin{equation}\label{pqcondition}
\fib^*\alpha\big(\underbrace{\xi_1~~,~~\dots~~,~~\xi_{p'}}_{\in\,\spn_\C(\bar \omega^2 \eta_- + \eta_+,\partial_{\omega}) },\underbrace{\xi_{p'+1}~~,~~\dots~~,~~\xi_{k}}_{\in\,\spn_\C(\omega^2\eta_++\eta_-,\partial_{\bar \omega})}\big)=0,\quad \text{for } p'\neq p.
\end{equation}
Using \eqref{taugamma}, the condition in the preceding display can be verified for the forms $\varepsilon^{p,q}_h$ defined in Table \ref{table} and this implies that $\fib_*\varepsilon_h^{p,q}$ is indeed a $(p,q)$-form. Consequently, the map $\Phi^{p,q}$ is well-defined as in \eqref{phipq}; it is clearly injective and we now demonstrate that it is onto. Let $\alpha\in \Omega^{p,q}(U)$, then the $\mathbf V$-invariant form $\fib^*\alpha$ can be written in terms of (wedge products of) the forms $\tau,\bar \tau,\gamma,\bar \gamma$. Using \eqref{pqcondition}, one sees that all but possibly $N= {2 \choose p}\cdot {2 \choose q}$ coefficients are zero, or more precisely, that
\begin{equation*}
\fib^*\alpha = \varepsilon_h^{p,q}\text{ for some } h=(h_1,\dots,h_N)\in C^\infty(\fib^{-1}(U))^N.
\end{equation*}
It remains to show that $h_1,\dots,h_N\in \mathcal I^{p-q}(U)$ and we do this exemplarily for the case $(p,q)=(1,2)$, where $N=2$ and $\fib^*\alpha = h_1 \bar \tau\wedge \tau \wedge \gamma + h_2 \bar \gamma \wedge \tau\wedge \gamma$. Then
\begin{equation}
0 = \mathcal L_\mathbf V(\fib^*\alpha)  = (\mathbf V h_1+ ih_1) \cdot \bar \tau\wedge \tau \wedge \gamma + (\mathbf V h_2+ ih_2) \cdot \bar \gamma \wedge \tau \wedge \gamma,
\end{equation}
which implies $h_1,h_2\in \ker (\V+i) = \mathcal I^{-1}(U).$

To verify the intertwining property, let $\alpha = \Phi^{p,q}(h)\in \Omega^{p,q}(U)$. Then the $(p,q+1)$-form $\fib^*(\bar \partial \alpha)$
is uniquely determined by its action on  vector fields $\xi_1,\dots,\xi_{k+1}$ with $\bar \xi_1,\dots,\bar \xi_{p'},\xi_{p'+1},\dots,\xi_{k+1}\in \spn_\C(\omega^2\eta_++\eta_-,\partial_{\bar \omega})$ through
\begin{equation}
\fib^*(\bar \partial \alpha)(\xi_1,\dots,\xi_{k+1}) = \begin{cases}
	\fib^*(d\alpha)(\xi_1,\dots,\xi_{k+1})&p'=p;\\
	0 & p'\neq p.
\end{cases}
\end{equation}
We claim that $\delta_h^{p,q}(\xi_1,\dots,\xi_{k+1})$ also satisfies this property; for $p'\neq p$ this is clear and for $p'=p$ this amounts to showing
\begin{equation}
\delta_{h}^{p,q}(\xi_1,\dots,\xi_{k+1}) = d\epsilon_h^{p,q}(\xi_1,\dots,\xi_{k+1}), 
\end{equation}
which can be done case by case, using the commutator formulas \eqref{equation:commutateur2} to compute the differential on the right---see also the proof of \cite[Lemma 4.2]{Bohr-Paternain-21}, where this is done in detail for $p=0$. 
\end{proof}

We turn this proposition into a definition in order to obtain a Dolbeaut complex for $(p,q)$-forms that are {\it smooth up to the boundary} of $Z$. That is, if $U\subset Z$ is an open set with $U\cap \partial Z\neq \emptyset$, then we simply {\it define}
\begin{equation}\label{pquptoboundary}\Omega^{p,q}(U):=\mathcal I^{p-q}(U)^{{2\choose p}{2\choose q}},\quad \bar \partial := D^{p,q}.\end{equation}
By abuse of notation, we keep writing $(p,q)$-forms as push-forwards. For example, a $(0,1)$-form that is smooth up the boundary will be written as
\begin{equation}\label{example01}
\fib_*(h_1 \tau + h_2 \gamma) \in \Omega^{0,1}(U),\quad h_1,h_2\in \mathcal I^{-1}(U).
\end{equation} 
However, the reader is cautioned that for $U\cap \partial Z\neq \emptyset$ this expression does {\it not} make sense as a differential form in $\Omega^{p+q}(U)$.

 \begin{remark} In \cite{Bohr-Paternain-21} the isomorphism in the preceding proposition was understood implicitly. For example, the $(0,1)$-form in \eqref{example01} was written simply as tuple $(h_1,h_2)$. While the formal notation chosen here might be slightly cumbersome, it is convenient to distinguish forms of different bi-degrees.
 \end{remark}

 \begin{remark}
 If $L\rightarrow Z$ is a complex line bundle then the space of twisted $(p,q)$-forms $\Omega^{p,q}(U,L)$ can be defined as in \eqref{pquptoboundary} by replacing $\mathcal I^m(U)$ with
 \begin{equation}
 \mathcal I^m(U,L)=\mathcal I^m(U)\otimes C^\infty(U,L).
 \end{equation}
 If $U\cap \partial Z =\emptyset$, then the same arguments as in Proposition \ref{intertwine} show that this agrees with the standard definition. If $U\cap \partial Z \neq \emptyset$, the elements of $\Omega^{p,q}(U,L)$ are thought of as twisted $(p,q)$-forms that are smooth up to the boundary, bearing in mind the same caveat as discussed below \eqref{example01}.

\end{remark}

We conclude this section with two lemmas that will be used below.

\begin{lemma}\label{smoothtrace}
Let $m\in \Z$ and $u\in C^\infty(SM)$ with $u_k=0$ for $k<m$. Then $h=\sum_{k\ge m}u_k \omega^k$ is a well defined function in $\mathcal I^m(Z)$ and the map $u\mapsto h$ fits into a split exact sequence:
\begin{equation*}
0\rightarrow \oplus_{k\ge m} C^\infty(M,\Omega_k) \rightarrow \mathcal I^m(Z) \xrightarrow{ \partial_{\bar \omega}} \mathcal I^{m-1}(Z)\rightarrow 0
\end{equation*}
\end{lemma}

This lemma appears as Lemma 4.4.4 in \cite{Bohr-thesis} and completes Proposition 4.4 in \cite{Bohr-Paternain-21} by proving also the stated mapping properties of $\partial_{\bar \omega}$; we include a proof of the latter for the reader's convenience.

\begin{proof}[Proof (Mapping properties of the  $ \partial_{\bar \omega}$-operator)] 
The fact that  $\partial_{\bar \omega}$ maps $\mathcal I^m(Z)$ into $\mathcal I^{m-1}(Z)$ follows from the observation that $\partial_{\bar \omega} \in \ker (\mathcal L_\mathbf V +i)$. To show surjectivity, let $f\in \mathcal I^{m-1}(Z)$ and extend it to a compactly supported, smooth function on $SM\times \C$. Then
$
g(x,v,\omega) = \frac{1}{2\pi i} \int_\C f(x,v,\zeta)(\omega-\zeta)^{-1} d\bar \zeta \wedge d\zeta
$ defines a solution $g\in C^\infty(SM\times \D)$ to $\partial_{\bar \omega}  g= f$. Using the spectral decomposition for $-i\mathbf V$ (similar to the one in  \eqref{equation:decomposition} for $-iV$), we may write $g=\sum_{k\in \Z} g_k$ with $g_k\in \mathcal I^k(Z)$. Arguing similarly as in the proof of Lemma 6.1.3 in \cite{Paternain-Salo-Uhlmann-book}, one shows that
\begin{equation*}
\partial_{\bar \omega} g_k = (\partial_{\bar \omega} g)_{k-1}\in \mathcal I^{k-1}(Z),\quad k\in \Z,
\end{equation*}
such that $h=g_m\in \mathcal I^m(Z)$ is the desired preimage of $f$ under $\partial_{\bar \omega}$. 
\end{proof}

\begin{lemma}[Lifted volume form]\label{lem:volume} Let $\zeta=\fib_*(\bar \tau \wedge \tau)\in \Omega^{1,1}(Z)$ and let $\iota_M\colon M\rightarrow Z$ be the inclusion as $0$-section. Then\begin{equation}
\iota_M^*\zeta = -2i~ d\mathrm{Vol}_g
\end{equation} 
\end{lemma}

\begin{proof}
Let $x\in M$ and consider an oriented orthonormal basis $\{v,v^\perp\}$ of $T_xM$. The horizontal lift of this basis to $T_{(x,v,0)}(SM\times \D)$ is $\{X(x,v),H(x,v)\}$ and thus
\begin{equation}
(\iota_M^*\zeta)_x(v,v^\perp) = \bar \tau \wedge \tau(X,H)\vert_{(x,v,0)}= -2i \cdot \bar \tau \wedge \tau(\eta_+,\eta_-)\vert_{(x,v,0)} = -2i,
\end{equation}
where we first used that $\eta_\pm = \frac 12 (X\mp i H)$ and then applied \eqref{taugamma} to compute $\tau(\eta_-)=1$ and $\tau(\eta_+)=0$ for $\omega=0$.
\end{proof}

\section{Invariant distributions as boundary values of holomorphic functions on twistor space}

In this section we prove Theorems \ref{mainthm1} and \ref{mainthm2}, starting with an underlying regularity result that is of independent interest. The proofs are based on microlocal analysis and we refer the reader to Appendix \ref{appendix} for a quick introduction.

\label{section:invariant}

\subsection{A regularity statement} 


We first prove the following statement on the regularity of invariant distributions. 

\begin{theorem}
\label{theorem:regularity}
Let $u \in \mc{D}'(SM)$ such that $(X+\pi^*\phi) u \in C^\infty(SM)$, where $\phi \in C^\infty(M)$. Then, the following assertions are equivalent:
\begin{enumerate}[label=\emph{(\roman*)}]
\item There exists an open conic neighbourhood $U$ of $\HH^*$ such that $\WF(u) \cap U = \emptyset$,
\item For all $k \in \Z$, $u_k \in C^\infty(M,\Omega_k)$,
\item $u_0 $ and $u_1$ are smooth.
\end{enumerate}
Moreover, the neighbourhood $U$ can be chosen independent of $u$ (it only depends on the metric $g$).
\end{theorem}

Observe that Theorem \ref{theorem:regularity} does not rely on any curvature assumption. The only key property on which it relies is the twist property of the vertical bundle, see \S\ref{sssection:twist}.

\begin{proof}
It is clear that (ii) implies (iii). Assume (iii) holds and note that the equation $(X+\pi^*\phi)u=f$ gives the recurrence relation:
\[\eta_{-}u_{k+1}+\eta_{+}u_{k-1}+\phi\,u_k=f_{k},\;\;k\in\mathbb{Z}.\]
An induction argument using the ellipticity of the operators $\eta_\pm$ and smoothness of $f$ gives that $u_k$ is smooth for all $k\in \mathbb{Z}$ provided
that $u_0$ are $u_{1}$ are smooth. 
The fact that (i) implies (ii) is because ${\pi_k}_*$ only selects wavefront set in $E_0^* \oplus \HH^*$, that is, $u_k = \pi_k^*{\pi_k}_* u$ and
\[
\begin{split}
\WF({\pi_k}_*u)  \subset \big\{ (x,\xi) \in T^*M \setminus \left\{0\right\} ~|~ \exists v \in S_x M, ( (x,v), \underbrace{d\pi^{\top}\xi}_{\in \HH^*\oplus E_0^*}, \underbrace{0}_{\in \V^*}) \in \WF(u) \big\} = \emptyset,
\end{split}
\]
see \cite[Lemma 3.2.11]{Lefeuvre-book} for a proof of the previous inclusion. (Heuristically, this can easily understood as ${\pi_k}_* u$ consists in integrating over the $S^1$ fibers of $SM$ so it ``kills'' all singularities tangential to the fibers; only the singularities that are tangential to the base direction remain after integration.)

We now prove that (ii) implies (i). First of all, by standard elliptic regularity, the equation $(X+\pi^*\phi) u = f$ implies that $\WF(u) \subset \Sigma$. Then, using \eqref{equation:szego-commutation}, and $[\mc{S},\pi^*\phi]=0$, we obtain
\begin{equation}
\label{equation:decomp-szego}
\begin{split}
(X+\pi^*\phi) \mc{S} u  = \mc{S} (X+\pi^*\phi) u + [X,\mc{S}] u  = \mc{S} (X+\pi^*\phi) u + \eta_- u_0 - \eta_+ u_{-1}.
\end{split}
\end{equation}
Observe that both terms on the right-hand side are smooth: indeed, $\mc{S} (X+\pi^*\phi)u$ is smooth since $(X+\pi^*\phi)u$ is, and the second term is a differential operator of order $1$ acting on the smooth sections $u_0+u_{-1}$, so it is smooth. Moreover, the wavefront set $\WF(K_\mc{S})$ of the Schwartz kernel $K_{\mc{S}}$ of $\mc{S}$ can be computed (it is just the usual Szeg\H{o} projector on each fibre $S_xM \simeq \Ss^1$) and turns out to be equal to
\begin{equation}
\label{equation:wf-szego}
\WF(K_{\mc{S}}) = \left\{ (v,\xi,v,-\xi) ~|~ (v,\xi) \in T^*(SM) \setminus \left\{0\right\}, \langle \xi,V(v)\rangle \geq 0\right\},
\end{equation}
see \cite[Lemma 3.10]{Guillarmou-17-1} for a proof.
As a consequence, we get by \eqref{equation:wf-szego} and the usual composition rules for the wavefront set (see \cite[Theorem 8.2.13]{Hormander-90}) that
\begin{equation}
\label{equation:wf-szego2}
\WF(\mc{S} u) \subset \left\{(v,\xi) \in T^*(SM) \setminus \left\{0\right\}, \langle \xi,V(v)\rangle \geq 0\right\}.
\end{equation}

We let $A \in \Psi^0(SM)$ be a pseudodifferential operator with wavefront set on a small conic neighbourhood of a point $(x_0,v_0, \xi_0) \in \HH^* \setminus \left\{0\right\}$ and elliptic on a slightly smaller conic neighbourhood of $\HH^*$. By the twist property of the cohorizontal bundle (see \S\ref{sssection:twist}), there exists a time $t_0 \in \R$ such that $\Phi_{t_0}(\WF(A)) \subset \left\{\langle\xi, V\rangle < 0\right\}$. Taking $B \in \Psi^0(SM)$, microlocally equal to the identity on a conic neighbourhood of $\Phi_{t_0}(\WF(A))$ and with wavefront set on a slightly larger conic neighbourhood (note that all this can be performed so that $\WF(B) \subset \left\{\langle \xi,V\rangle < 0\right\}$), we obtain by Lemma \ref{lemma:propagation} for $P=\pm iX$ the existence of $B_0 \in \Psi^0(SM)$, with wavefront set contained in a conic neighbourhood of $\cup_{t \in [0,t_0]} \Phi_t(\WF(A))$, such that for all $s \in \R, N > 0$, there exists a constant $C >0$ such that for all $f \in C^\infty(SM)$,
\begin{equation}
\label{equation:bound}
\|A f\|_{H^s} \leq C \left( \|B f\|_{H^s} + \|B_0 (X+\pi^*\phi)f\|_{H^s} + \|f\|_{H^{-N}} \right).
\end{equation}
Moreover, the left-hand side is defined and bounded by \eqref{equation:bound} whenever all terms on the right-hand side are finite and $f \in \mc{D}'(SM)$. In particular, we can apply \eqref{equation:bound} to $f = \mc{S}u$: by construction, the right-hand side is indeed finite since $\mc{S} u$ is microlocally smooth on $\WF(B)$ and $(X+\pi^*\phi)\mc{S} u$ is a smooth function. Hence, we obtain:
\begin{equation}
\label{equation:bound2}
\|A \mc{S} u\|_{H^s} \leq C \left( \|B \mc{S} u\|_{H^s} + \|B_0 (X+\pi^*\phi) \mc{S} u\|_{H^s} + \|\mc{S}  u\|_{H^{-N}} \right)
\end{equation}
Since \eqref{equation:bound2} can be applied with an arbitrary $s \in \R$, we deduce that $\mc{S} u$ is microlocally smooth on $\WF(A)$. Note that, since $B$ has wavefront set in $\left\{\langle\xi,V\rangle<0\right\}$, the operator $B \mc{S} \in \Psi^{-\infty}(SM)$ is smoothing and $\|B\mc{S} u\|_{H^s} \leq C\|u\|_{H^{-N}}$. Hence, we get from \eqref{equation:bound2} the new bound:
\begin{equation}
\label{equation:bound3}
\|A \mc{S} u\|_{H^s} \leq C \left(\|B_0 (X+\pi^*\phi) \mc{S} u\|_{H^s} + \|u\|_{H^{-N}} \right).
\end{equation}
Finally, it suffices to observe that a similar argument can be applied to $(\mathbbm{1}-\mc{S})u$, allowing to conclude that $u$ is microlocally smooth on $\WF(A)$. This ends the proof.
\end{proof}

The theorem has the following noteworthy corollary:

\begin{corollary}
\label{corollary:gabriel}
Let $u \in \mc{D}'(SM)$ be such that $X u =0$ and $u_0+u_{1} \in C^\infty(M,\Omega_{0}\oplus\Omega_{1})$. Then the singular support of $u$ is a closed invariant set in $SM$ free of conjugate points.\label{corollary:freeconjugate}
\end{corollary}

\begin{proof} By Theorem \ref{theorem:regularity} we know that $\WF(u)\cap \HH^* = \emptyset$. By propagation of singularities \cite[Theorem 26.1.1]{Hormander-4}, the equation $Xu=0$ gives that $\WF(u)\subset \Sigma$ and $\WF(u)$ is invariant under $(\Phi^{T^*(SM)}_{t})_{t\in\R}$. In particular, the singular support of $u$ is a closed set in $SM$ invariant under the geodesic flow $\varphi_t$. Moreover, if we take $v$ in the singular support, $\WF(u)$ gives rise to a Jacobi field along the geodesic determined by $v$ which never vanishes since  $\WF(u)\cap \HH^* = \emptyset$. By Sturm comparison it follows that the geodesic determined by $v$ is free of conjugate points.
\end{proof}

In particular, in the case of a positively-curved metric $g$ on $\Ss^2$, every geodesic has conjugate points, and thus Corollary \ref{corollary:gabriel} shows that a flow-invariant distribution $u$ with $u_0,u_1 \in C^\infty(M,\Omega_0\oplus\Omega_1)$ is actually smooth. This is a non-quantitative version of Theorem \ref{theorem:sphere-intro}.

\subsection{Unital algebra of invariant distributions}

Define the following subspace of invariant distributions:
\begin{equation}
\label{equation:laurent-distributions}
\begin{split}
\mc{L}'(SM) := \big\{ u \in \mc{D}'(SM) ~|~ Xu = 0, \exists k_0 \in \Z_{\geq 0}, u = \sum_{k \geq -k_0} u_k\big\}   \subset \mc{D}'(SM).
\end{split}
\end{equation}
Such distributions are exactly the ones admitting a holomorphic extension to the sliced twistor space $Z\backslash 0$, cf.\,\S\ref{section:twistor}. It is natural to equip $\mc{L}'(SM)$ with the topology inherited from $\mc{D}'(SM)$. We will prove that the following holds:

\begin{corollary}
\label{corollary:multiplication}
The wavefront set of a distribution $u\in \mc{L}'(SM)$ satisfies
\begin{equation}
\WF(u)\subset  \mathcal C\subset \{(v,\xi)\in \Sigma: \langle \xi,V(v)\rangle >0\}
\end{equation}
for a closed cone $\mathcal C$ that only depends on the metric. In particular,
$\mc{L}'(SM)$ is a unital algebra with continuous multiplication.
\end{corollary}

When there are no conjugate points, $\mc{C}$ is the closed cone between $E_s^*$ and $E_u^*$ in $\left\{ \langle \xi, V(v)\rangle  > 0 \right\}$, as in Figure \ref{figure:cones}. Note that, as a particular consequence of Corollary \ref{corollary:gabriel}, when $M=\Ss^2$ is equipped with a positively-curved metric, $\mc{L}'(SM) \subset C^\infty(SM)$ and the algebra property is therefore immediate. The proof of Corollary \ref{corollary:multiplication} is reminiscent of \cite[Corollary 3.11]{Guillarmou-17-1}.

\begin{proof}
Any $u\in \mc L'(SM)$ clearly satisfies $\WF(u)\subset \Sigma$. Using the ellipticity of $\eta_\pm$, we have $u_k \in C^\infty(M,\Omega_k)$ for all $k \in \Z$, which implies that $u -\mc{S} u = \sum_{k=-k_0}^{-1} u_k \in C^\infty(SM)$ and thus $\WF(u) \subset \left\{\langle\xi,V\rangle\geq 0\right\} \cap \Sigma$ by \eqref{equation:wf-szego2}. Finally,  Theorem \ref{theorem:regularity}(i) implies that $\WF(u)$ is contained in a closed cone $\mathcal C\subset \left\{\langle\xi,V\rangle> 0\right\} \cap \Sigma$, which is independent of $u$.
To see that $\mathcal L'(SM)$ is an algebra, notice first there is no $(v,\xi)\in \mathcal C$ with $(v,-\xi)\in \mathcal C$, such that multiplication of elements in $\mathcal L'(SM)$ is well-defined and continuous by \cite[Theorem 8.2.10]{Hormander-90}. It is a unital algebra, as it contains the constant function $\mathbf{1}$. This proves the corollary.
\end{proof}

\subsection{The trace map} 
We wish to consider fibrewise holomorphic invariant distributions on $SM$ as traces of holomorphic functions in the interior $Z^\circ$ of twistor space.
If $f\in C^\infty(Z^\circ)$, then $f^{r}(x,v)=f(x,rv)$ ($0\le r<1$) defines a smooth function on $SM$  and we define the trace of $f$ as 
\begin{equation}\label{ftrace}
f\vert_{SM} := \lim_{r\rightarrow 1^{-}} f^{r} \in \mathcal D'(SM),
\end{equation}
provided this limit exists in the sense of distributions. This limit is the content of the following proposition, which is stated for the spaces $\mathcal I^m(Z^\circ)$ ($m\in \Z$)  from \S\ref{section:dolbeaut}. In order to apply the proposition to a holomorphic function $f\in C^\infty(Z^\circ)$, one considers  $h=
\fib^*f\in \mathcal I^0(Z^\circ)$. Then $f^{r}=h\vert_{\omega=r}\in C^\infty(SM)$ and
 the limit in  \eqref{ftrace} is a special case of \eqref{dtrace} below. 
 
Below the space $C^N(M,\Omega_k)$ of $N$-times continuously differentiable sections of $\Omega_k$ is normed in a $k$-independent way by considering it as subspace of $C^N(SM)$ and requiring that the embedding $\pi_k^*\colon C^N(M,\Omega_k)\rightarrow C^N(SM)$  from \eqref{equation:identification} is isometric for all $k\in \Z,N\in \mathbb N$.

\begin{proposition}\label{prop_trace1} Let $m\in \N$ and suppose $h\in \mathcal I^m(Z^\circ)$ satisfies $\partial_{\bar \omega} h=0$. 
\begin{enumerate}[label=\emph{(\roman*)}]
\item Then $h = \sum_{k\ge m} u_k\omega^{k-m}$ for coefficients $u_k\in C^\infty(M,\Omega_k)$ and this series converges in the topology of $C^\infty(SM\times \D^\circ)$; 
\smallskip
\item for all $N\in \mathbb N$ the following polynomial growth conditions are equivalent:
\smallskip
\begin{eqnarray}
\Vert u_k \Vert_{C^N(M)} &\le & C \langle k \rangle^p \,\quad\qquad \text{ for some } C,p>0;\label{growth1}\\
\Vert h\vert_{\omega=r} \Vert_{C^N(SM)} &\le& C(1-r)^{-p} \quad \text{ for some } C,p>0 \label{growth2}.
\end{eqnarray}
\item if these growth conditions are satisfied for some $N\in \mathbb N$, then the following  limits exist in $\mathcal D'(SM)$ and they agree:
\begin{equation}
\label{dtrace}
 \lim_{r\rightarrow 1^-} h\vert_{\omega=r} =\lim_{P\rightarrow \infty} \sum_{k=m}^P u_k  \in \mathcal D'(SM);
\end{equation}
\item if this limit is zero, then $h$ must vanish identically.
\end{enumerate}
\end{proposition}

We first prove the following lemma:

\begin{lemma}\label{smoothdecay} For every $N,q\in \mathbb N$ there exists a constant $C>0$ such that \begin{equation*}
\Vert \varphi_k \Vert_{C^N(M)} \le C \langle k\rangle^{-q} \Vert \varphi \Vert_{C^{N+q}(SM)}
\end{equation*}
for all $\varphi\in C^\infty(SM)$ and all $k\in \Z$. In particular, $\Vert \varphi_k\Vert_{C^N{(M)}}=O(\langle k\rangle^{-\infty})$.
\end{lemma}

\begin{proof}
We have $\vert k \vert^q \Vert \varphi_k \Vert_{C^N(M)} =  \Vert (ik)^q \varphi_k \Vert_{C^N(M)}=\Vert V^q \varphi_k \Vert_{C^N(SM)} \lesssim \Vert \varphi_k \Vert_{C^{N+q}(SM)}$ and this yields a uniform bound as in the lemma, at least for all $k\neq 0$, where $\langle k \rangle \le \sqrt 2 \vert k \vert$. The estimate for $k=0$ is trivial and does not affect uniformity.
\end{proof}

\begin{proof}[Proof of Proposition \ref{prop_trace1}]
For (i) let $h\in \mathcal I^m(Z^\circ)\cap \ker \partial_{\bar \omega}$; then for fixed $(x,v)\in SM$ the function $h(x,v,\cdot)$ is holomorphic on $\{\vert \omega \vert<1\}$ and thus it has a series expansion $
h(x,v,\omega)=\sum_{k\ge m} u_k(x,v) \omega^{k-m}
$ with coefficients  $u_k(x,v)\in \C$ given in terms of the Cauchy integral formula, that is,
\begin{equation}\label{cauchyintegral}
 u_k(x,v) = \frac{1}{2\pi i} \int_{\vert \zeta \vert = r}\frac{h(x,v,\zeta)}{\zeta^{k-m+1}} d\zeta,\quad 0<r<1.
\end{equation}
This shows immediately that $u_k\in C^\infty(SM)$. It is clear that the power series for $h$ converges pointwise; we now show that the partial sums and their derivatives are Cauchy in $L^\infty(SM\times \{\vert \omega \vert\le s\})$ for all $0<s<1$. Choose $r\in (s,1)$, then as $P,P'\rightarrow \infty$,
\begin{equation}
\sum_{k=P+1}^{P'} u_k \omega^{k-m} = \frac {1}{2\pi i}  \sum_{k=P+1}^{P'} (s/r)^k \int_0^{2\pi} \frac{s^{-k}\omega^{k-m}}{e^{i(k-m+1)\theta}} h(x,v,re^{i\theta}) r^{m}  d\theta \rightarrow 0,
\end{equation}
uniformly for all $(x,v,\omega)\in SM\times \D$ with $\vert \omega \vert \le s$. As $h$ is assumed to be smooth, the same holds true after taking an arbitrary number of derivatives and thus the power series indeed converges in the topology of $C^\infty(SM\times \D^\circ)$.
We may thus compute the derivative $\mathbf{V} h$ term by term and observe that
\begin{equation}\label{huseries1}
0 = (\mathbf{V} - i m)h = \sum_{k\ge m} (V-im)u_k  \omega^{k-m} - i(k-m) u_k \omega^{k-m},
\end{equation}
which implies that $(V-ik)u_k=0$ and thus $u_k\in C^\infty(M,\Omega_k)$, as desired.

For (ii), let us first assume that \eqref{growth1} holds true. Then for $1/2<r<1$,
\begin{equation}
\Vert h\vert_{\omega =r} \Vert_{C^N(SM)} \le \sum_{k\ge m} r^{k-m} \Vert u_k\Vert_{C^N(M)}  \lesssim  \sum_{k\ge m} r^{k-m} \langle k\rangle^{p} \lesssim 1 + \sum_{k\ge 0} r^{k} k^p,
\end{equation}
with constants that are uniform in $r$. But $r^k k^p \le r^k \frac{(k+p)!}{p!} = (d/dr)^pr^{k+p}$, such that the previous display can be continued with 
\begin{equation}
\lesssim 1 + (d/dr)^p\left[r^p(1-r)^{-1}\right] \lesssim (1-r)^{-(p+1)},
\end{equation}
which yields \eqref{growth2}. Conversely, \eqref{growth2} and  \eqref{cauchyintegral} imply that \begin{equation}
\Vert u_k \Vert_{C^N(M)} \le r^{-k+m} \Vert h\vert_{\omega=r}\Vert_{C^N(SM)}\lesssim r^{-k+m} (1-r)^{-p},
\end{equation}
for some $p\in \N$, with constants that are uniform in $k\in \N$ and $0<r<1$. Choosing $r=1-1/k$, we may continue the estimate with
\begin{equation}
\lesssim (1-1/k)^{-k+m}
 k^p \lesssim \langle k \rangle^p,
\end{equation}
where we used that $(1-1/k)^{-k+m}$ has a limit as $k\rightarrow \infty$, such that it is bounded. 

For (iii) consider a test function $\varphi\in C^\infty(SM)$. Then by \eqref{growth1} and Lemma \ref{smoothdecay} (for $N=0$),  there exists $p\in \N$ such that for all $q\in \N$,
\begin{equation}
\vert\langle u_k ,\varphi_k\rangle\vert \lesssim \Vert u_k \Vert_{L^\infty(SM)} \Vert \varphi_k \Vert_{L^\infty(SM)} \lesssim   \langle k \rangle^p \cdot \langle k\rangle^{-q} \Vert \varphi\Vert_{C^q(SM)}.
\end{equation}
Choosing $q$ large enough, this implies that the series $\sum_{k\ge m} \langle u_k,\varphi_k\rangle$ converges absolutely, with limit satisfying
$
\vert \sum_{k\ge m} \langle u_k,\varphi_k\rangle \vert \lesssim \Vert \varphi \Vert_{C^q(SM)}.
$
Hence, $\varphi\mapsto \sum_{k\ge m} \langle u_k,\varphi_k\rangle$ defines a distribution $u\in \mathcal D'(SM)$. By the dominated convergence theorem,
\begin{equation*}
\langle h\vert_{\omega=r} ,\varphi\rangle = \sum_{k\ge m} r^k \langle u_k,\varphi_k\rangle \xrightarrow{r\rightarrow 1} \langle u,\varphi\rangle \quad \text{and} \quad \langle \sum_{k=m}^P u_k,\varphi\rangle =  \sum_{k=m}^P \langle u_k,\varphi_k\rangle \xrightarrow{P\rightarrow \infty} \langle u,\varphi\rangle,
\end{equation*}
thus proving (iii). Finally, if $u=0$ then  $u_k=0$ for all $k\ge m$ and $h\equiv 0$ follows from the expansion in (i). This concludes the proof.
\end{proof}

\begin{corollary} \label{traceiso1} The trace map from Proposition \ref{prop_trace1} is a linear isomorphism
\begin{equation*}
\begin{split}
&\left \{h\in \mathcal I^m(Z^\circ):
\begin{array}{l}
 \partial_{\bar \omega}h = 0 \text{ and } 
\eqref{growth2} \\ \text{holds for all } N\in \mathbb N
\end{array}
\right\} \\
& \hspace{4.5cm} \longrightarrow
\left \{u\in \mathcal D'(SM):
\begin{array}{l}
 u_k = 0 \text{ for } k<m,  \text{ its Fourier}\\
 \text{modes are smooth and satisfy}
\\
\text{\eqref{growth1} for all } N\in \mathbb N
\end{array}
\right\}.
\end{split}
\end{equation*}
\end{corollary}

\begin{proof}
	If $u = \lim_{r\rightarrow 1^{-}} h\vert _{\omega=r}$ is the trace of some $h\in \mathcal I^m(Z^\circ)\cap \ker \partial_{\bar \omega}$, then its nonzero Fourier modes $u_k$ $(k\ge m)$ arise as in Proposition \ref{prop_trace1}(i) and in particular they are smooth and satisfy \eqref{growth1}. The injectivity of the trace map follows from part (iv) of the proposition.
	Vice versa, suppose that $u\in \mathcal D'(SM)$ lies in the space on the right. We claim that \[
h= \sum_{k\ge m}u_k \omega^{k-m} 
\] 
converges in the topology of $C^\infty(SM\times \D^\circ)$. Pointwise convergence of this series follows immediately from \eqref{growth1} (for $N=1$), so it remains to show that the partial sums form a Cauchy sequence in $C^K(SM\times \{\vert \omega\vert \le r\})$ for all $0<r<1, K\in \mathbb N$. As the partial sums are holomorphic in $\omega$, their $C^K$-norms can be expressed in terms of differential operators of the form $Q\partial_{\omega}^q$, where $Q$ is a differential operator on $SM$ and $\mathrm{ord}\,Q + q \le K$. Let $p\in \N$ be the exponent from \eqref{growth1} for $N=\mathrm{ord}\,Q$, then for $P'\ge P\ge q+m$ we have
\[
\begin{split}
\Big \Vert Q \partial_\omega^q& \sum_{k=P+1}^{P'} u_k \omega^{k-m}  \Big\Vert_{L^\infty(SM\times \{\vert \omega \vert \le r\})} 
 \\
&= \frac{(k-m)!}{(k-m-q)!}\Big\Vert  \sum_{k=P+1}^{P'} (Qu_k) \omega^{k-m-q}  \Big\Vert_{L^\infty(SM\times \{\vert \omega \vert \le r\})}\\
&\lesssim
 \sum_{k=P+1}^{P'} \Vert u_k \Vert_{C^N(M)} \cdot \Vert \omega^{k-m-q} \Vert_{L^\infty(\{\vert \omega \vert \le r\})} 
 \lesssim  \sum_{k=P+1}^{P'} \langle k \rangle^p \cdot r^{k-m-q} \xrightarrow{P,P'\rightarrow \infty} 0,
\end{split}
\]
and this proves the claim. By the same computation as in \eqref{huseries1} one shows that $\mathbf{V} h =i m h$, such that $h\in \mathcal I^m(Z^\circ)\cap \ker \partial_{\bar \omega}$. Finally, by part (ii) of Proposition \ref{prop_trace1}, property \eqref{growth2} is satisfied for all $N\in \mathbb N$.
\end{proof}

\begin{proposition} \label{wftest}For $u \in \mathcal D'(SM)$ consider the following three properties:
\begin{enumerate}[label=\emph{(\roman*)}]
\item $\WF(u)\cap (E_0^*\oplus \mathbb H^*) = \emptyset$;
\smallskip
\item $\pi_*(u\varphi)\in C^\infty(M)$ for all $\varphi\in C^\infty(SM)$;
\smallskip
\item the Fourier modes of $u$ are smooth and satisfy \eqref{growth1} for all $N\in \mathbb N$.
\end{enumerate}
Then
$
\text{(i)} \Longrightarrow \text{(ii)}  \Longleftrightarrow \text{(iii)}. 
$
\end{proposition}

\begin{proof}
If (i) holds true, then also $\WF(u\varphi)\cap (E_0^*\oplus \mathbb H^*)=\emptyset$ and thus
$\WF(\pi_*(u\varphi)) = \emptyset$ by the same argument as in the proof of Theorem \ref{theorem:regularity}; this implies (ii).  For the equivalence between (ii) and (iii) we may assume without loss of generality that $\supp u\subset SM_1^\circ$ for a compact embedded disk $M_1\subset M$. Then the following linear map is continuous:
\begin{equation*}
\varphi \mapsto \pi_*(u\varphi),\quad C^\infty(SM_1)\rightarrow \mathcal D'(M_1^\circ).
\end{equation*}
Define $\psi^{k}\in C^\infty(M_1)$ in isothermal coordinates as $\psi^k(x,\theta)=e^{-ik\theta}$, then the $k$th Fourier mode of $u$ can be identified with the compactly supported distribution
$
u_k = \pi_*(u\psi^k)
$. 

Assuming (ii), it immediately follows that the Fourier modes of $u$ are smooth. Further, the map in the preceding display then takes values in $C^\infty(M_1)$ and, using the closed graph theorem for Fr{\'e}chet spaces, it must also map continuously into $C^\infty(M_1)$. This means that
\begin{equation*}
\forall N\in \N~\exists p\in \N:\quad \Vert u_k\Vert_{C^N(M)} =\Vert \pi_*(u\psi^k)\Vert_{C^N(M_1)} \lesssim_N \Vert \psi^k\Vert_{C^{p}(SM_1)} \lesssim \langle k \rangle^{p},
\end{equation*}
with implicit constant independent of $k\in \Z$. Hence (iii) holds true. 

Conversely, assume that (iii) holds true and $\varphi\in C^\infty(SM_1)$. Then formally we have
\begin{equation*}
\pi_*(u\varphi) = \pi_*\left(u \sum_{k\in \Z} \varphi_k \psi^k\right)  = \sum_{k\in Z}  u_k \varphi_k
\end{equation*}
and in order to prove (ii) we need to show that the series converges in $C^\infty(M_1)$. Let $N\in \mathbb N$, then using (iii) as well as Lemma \ref{smoothdecay}, there exists $p\in \N$ such that
\begin{equation*}
\Vert u_k\varphi_k \Vert_{C^N(M_1)} \le \Vert u_k\Vert_{C^N(M)} \Vert \varphi_k \Vert_{C^N(M_1)} \lesssim_{N,q} \langle k\rangle^{p}\langle k\rangle^{-q} 
\end{equation*}
for all $q\in \N$. Choosing $q\ge p+2$ gives the desired convergence.
\end{proof}

\subsection{The algebra isomorphism} We can now complete the proof of Theorem \ref{mainthm1}, which can be rephrased as saying that the trace map from \eqref{ftrace} is a continuous algebra monomorphism
\begin{equation}
\mathcal A_\pol(Z)\rightarrow \mathcal L'(SM),\quad f\mapsto f\vert_{SM},
\end{equation}
with image $\mathcal L'_0(SM)=\{u\in \mathcal L'(SM): u \text{ is fibrewise holomorphic} \}$.

\begin{proof}[Proof of Theorem \ref{mainthm1}]
Let $f\in \mathcal A_\pol(Z)$ and $h=\fib^*f\in \mathcal I^0(Z^\circ)$. Expanding $h = \sum_{k\ge 0}u_k\omega^k$ as in Proposition \ref{prop_trace1}(i), holomorphicity of $f$ implies that
\begin{equation}
0 = (\omega^2 \eta_+ + \eta_-)\sum_{k\ge 0} u_k \omega^k = \sum_{k\ge 0} (\eta_+u_{k-2} + \eta_-u_k) \omega^k,
\end{equation}
where we set $u_k=0$ for $k<0$. It follows that $\eta_+ u_{k-2} + \eta_-u_k = 0$ for all $k\in \Z$ and as the $u_k$'s are the Fourier modes of $u=f\vert_{SM}$, we see that $Xu=0$. This shows that the $f\mapsto f\vert_{SM}$ indeed maps into $\mathcal L_0'(SM)$.
The fact that the trace map is an algebra homomorphism follows from the description \eqref{dtrace} and injectivity follows from Proposition \ref{prop_trace1}(iv).

We now show that any $u\in \mathcal L_0'(SM)$ lies in the image of the trace map. By Corollary \ref{corollary:multiplication}, the wavefront set of $u$ satisfies $\WF(u) \cap (E_0^*\oplus\mathbb H^*)=\emptyset$; hence Proposition \ref{wftest} and Corollary \ref{traceiso1} together imply that $u=f\vert_{SM}$ for some $f\in  C^\infty(Z^\circ)$ with $h=\fib^*f\in \mathcal I^0(Z^\circ) \cap \ker \partial_{\bar \omega}$ satisfying \eqref{growth2} for all $N\in \mathbb N$. The same computation as in the preceding display shows that $f$ is indeed holomorphic.
\end{proof}

\begin{remark}
The proof shows that functions $f\in \mathcal A_\pol(Z)$ automatically satisfy a stronger growth condition: if $f^r(x,v)=f(x,rv)$, then $\Vert f^r \Vert_{C^N(SM)}=O\left((1-r)^{-p}\right)$ for some $p\in \N$. For $N=0$ this is part of the definition of $\mathcal A_\pol(Z)$, but the case $N\ge 1$ follows a posteriori.
\end{remark}


%

%
%
%
%

\subsection{Some algebras of holomorphic functions on twistor space} We now proceed with the proof of Theorem \ref{mainthm2}. We start by computing the space of holomorphic functions on twistor space when $M=\Ss^2$.

\begin{proposition}
\label{proposition:s2}
Let $(\Ss^2,g)$ be the $2$-sphere equipped with an arbitrary metric $g$. Then:
\[
\mc{A}(Z_{\Ss^2}) = \mc{A}_{\mathrm{pol}}(Z_{\Ss^2}) = \mc{A}(Z_{\Ss^2}^\circ) = \C.
\]
\end{proposition}

\begin{proof}
Let $f \in \mc{A}(Z_{\Ss^2}^\circ)$. Writing $\fib^*f = \sum_{k \geq 0} u_k \omega^k \in C^\infty(SM \times \D^\circ)$ as in Proposition \ref{prop_trace1}, the equation $(\omega^2 \eta_++\eta_-)\fib^*f = 0$ reads $\eta_- u_0 = \eta_- u_1 = 0$ and $\eta_+ u_k + \eta_- u_{k+2}=0$ for all $k \geq 0$. Applying Lemma \ref{lemma:injectivity-eta}, we then easily get that $u_0 = c$ is constant, $u_k = 0$ for all $k \geq 1$. This proves the result.
\end{proof}

We now consider of the case of the flat torus.

\begin{proposition}
\label{proposition:t2}
Let $(\T^2,g_{\mathrm{flat}})$ be a flat $2$-torus. Then 
\[
\mathcal A(Z) \simeq \mathcal A(\D),\quad \mathcal A_\pol(Z) \simeq \mathcal A_\pol(\D),\quad \mathcal A(Z^\circ) \simeq \mathcal A(\D^\circ).
\]
Hence, $\mc{A}(Z_{\T^2,g_{\mathrm{flat}}}^\circ)$ does not separate points on $Z$ but $Z^\circ$ is holomorphically convex.
\end{proposition}


\begin{proof}
We identify $S\T^2 \simeq \T^2 \times \Ss^1$ and use the coordinates $(x, e^{i\theta})$. One can check by hand that for all $k \in \Z$, $\ker \eta_\pm|_{C^\infty(\T^2,\Omega_k)} = \C\cdot e^{i k \theta}$. Hence, given $f \in \mc{A}(Z_{\T^2,g_{\mathrm{flat}}}^\circ)$, the equation $(\omega^2 \eta_++\eta_-)\fib^*f = 0$ implies that $u_k(x,v) = c_k e^{i k \theta}$ for all $k \geq 0$, where $c_k \in \C$ is constant. As a consequence, we can map $f$ to the holomorphic function on $\D^\circ$ given by $\mu \mapsto \sum_{k \geq 0} c_k \mu^k$ and this mapping is an isomorphism. This shows that $\mathcal A(Z^\circ) \simeq \mathcal A(\D^\circ)$ and similarly, we obtain the other isomorphisms. However, observe that $\mc{A}(Z_{\T^2,g_{\mathrm{flat}}}^\circ)$ only separates point in the fibres $\D$ of the twistor space $Z \simeq \T^2 \times \D$ but it does not separate points in the base direction. Finally, $Z^\circ$ is holomorphically convex as $\D^\circ$ is.
\end{proof}

Eventually, we consider the case of Anosov surfaces. Let
\[
\mc{H}_m \simeq C^\infty(M,\Omega_m) \cap \ker \eta_-.
\]
be the space of holomorphic $m$-differentials (see \S\ref{section:fourier}, where this isomorphism is explained).
We will denote by $\mc{A}_\pol^m(Z)$ the space of holomorphic functions in $\mc{A}_\pol(Z)$ with Fourier degree $\geq m$. Recall that a Riemann surface $M$ is said to be \emph{hyperelliptic} if there exists a branched double cover $M \to \C P^1$. The following holds:

\begin{proposition}
Let $(M,g)$ be an Anosov surface. Then $\mc{A}(Z) = \C$ and both $\mc{A}_{\mathrm{pol}}(Z), \mc{A}(Z^\circ)$ are infinite-dimensional. Assuming moreover that the underlying Riemann surface $(M,[g])$ is non hyperelliptic, the spaces $\mathcal A_\pol^m(Z)$ ($m\ge 0$) fit into a short exact sequence
\begin{equation}
\label{equation:sequence}
0\rightarrow \mathcal A_\pol^{m+1}(Z) \rightarrow \mathcal A_\pol^{m}(Z) \rightarrow \mathcal H_m \rightarrow 0,\quad m\ge 0,
\end{equation}
where $\mathcal H_m$ is the space of holomorphic $m$-differentials of $(M,[g])$ (The condition on hyperellipticity can be dropped in the case $m=0$ and $m=1$). 
\end{proposition}

The last arrow in \eqref{equation:sequence} is simply given by taking the $m$-th Fourier mode $u_m \in C^\infty(M,\Omega_m)$ of the trace $f\vert_{SM}$ of a holomorphic function $u=f \in \mc{A}_\pol^m(Z)$. The fact that $X u = 0$ translates into $\eta_- u_m = 0$, that is, $u_m \in \mc{H}_m$.

\begin{proof}
That $\mc{A}(Z) = \C$ is immediate as any smooth flow-invariant function $f \in C^\infty(SM)$ is constant by ergodicity of the geodesic flow (more generally, $\mc{A}(Z) = \C$ as long as the geodesic flow is ergodic). The fact that the image of $\mathcal A_\pol^{m+1}(Z)$ in $\mathcal A_\pol^{m}(Z)$ is given by the kernel of the second arrow in \eqref{equation:sequence} is tautological. All that remains to be proven is that, under the assumptions, for any $a_m \in \mc{H}_m$, one can find $f \in \mathcal A_\pol^{m}(Z)$ with $m$-th Fourier mode given by $a_m$.

For $m=0$, $\mc{H}_0 = \C \cdot \mathbf{1}$ (the constant function). Hence, taking the function $f := \mathbf{1}$, we get a holomorphic function on $Z$ whose zeroth Fourier mode given by $\mathbf{1}$. Given $a_1 \in \mc{H}_1$, it was proved in \cite[Theorem 5.5]{Paternain-Salo-Uhlmann-14-2} that there exists a fibrewise holomorphic distribution $w \in \mc{D}'(SM)$ such that $Xw = 0$ and $w_1 = a_1$. Hence, taking $f \in \mc{A}^1_\pol(Z)$ such that $f|_{SM} = w$, we get the result. Note that, from there, we immediately get that $\mc{A}_\pol(Z)$ (hence $\mc{A}(Z^\circ)$) is always infinite-dimensional since it contains the infinite-dimensional space $\mathrm{span}(f^m, m \geq 0)$. Eventually, for $m \geq 2$, given $a_m \in \mc{H}_m$, using the non-hyperellipticity of $(M,[g]$), it was proved in \cite[Proof of Proposition 3.13]{Guillarmou-17-1} that there exists $w \in \mc{D}'(SM)$ such that $Xw=0$ and $w_m = a_m$. Taking as before $f \in \mc{A}^m_\pol(Z)$ with trace given by $w$ completes the proof.
\end{proof}

\section{Intermezzo: regularity on convex surfaces}

\label{section:intermezzo}

We now use the ideas of the previous paragraph to prove Theorem \ref{theorem:sphere-intro} on the regularity of invariant distributions on convex surfaces. In this section, it is assumed that $(M,g)$ is diffeomorphic to $\mathbb S^2$ and has positive curvature.
We will actually prove the following more general result:

\begin{theorem}
\label{theorem:sphere}
Let $\phi \in C^\infty(M)$. For all $s \in \R$, there exists a constant $C > 0$ such that the following inequality holds: for all $u \in C^\infty(SM)$,
\begin{equation}
\label{equation:stabilite}
\|u\|_{H^s} \leq C \left( \|(X+\pi^*\phi)u\|_{H^s} + \|u_0\|_{H^{s+1}}+\|u_1\|_{H^{s+1}}\right).
\end{equation}
Moreover, if $u \in \mc{D}'(SM)$ is a distribution, then for all $s \in \R$,
\[
(X+\pi^*\phi)u \in H^s(SM), u_0,u_1 \in H^{s+1}(SM) \implies u \in H^s(SM),
\]
and \eqref{equation:stabilite} holds. In particular, if $u \in \mc{D}'(SM)$, $(X+\pi^*\phi)u=0$ and $u_0+u_1 = 0$, then $u=0$.
\end{theorem}

Note that if $g$ is a Zoll metric on $M=\Ss^2$ (with period $2\pi$), then we have for all $\lambda \in \R\setminus \Z$:
\[
(e^{2i\pi\lambda}-1)^{-1} \int_0^{2\pi} e^{t(X+i\lambda)} |dt| (X+i\lambda) = \mathbbm{1}.
\]

As a consequence, we get for all $s \in \R$, that there exists a constant $C > 0$ such that for all $u \in C^\infty(SM)$:
\begin{equation}
\label{equation:bound-lambda}
\|u\|_{H^s} \leq C \|(X+i\lambda)u\|_{H^s}.
\end{equation}
This should be compared with \eqref{equation:stabilite} applied with $\phi := i\lambda$. Note that \eqref{equation:bound-lambda} generalises immediately to all $\phi \in C^\infty(\Ss^2)$ such that $\int_\gamma \phi \notin 2i\pi\Z$ for all great circles $\gamma \subset \Ss^2$.

\begin{proof}[Proof of Theorem \ref{theorem:sphere}]
The proof is divided in three steps. \\

\emph{Step 1. Injectivity statement.} We start by proving the injectivity statement which generalises \cite{Flaminio-92}:

\begin{lemma}
\label{lemma:inj}
Let $\phi \in C^\infty(M)$. Assume that $u \in \mc{D}'(SM)$ satisfies $(X+\pi^*\phi)u = 0$, $u_0 + u_1 = 0$. Then $u \equiv 0$.
\end{lemma}

\begin{proof}
Decomposing the equation $(X+\pi^*\phi)u = 0$ into each Fourier mode, we get for all $k \in \Z$:
\begin{equation}
\label{equation:decomp}
\eta_+ u_{k-1} + \phi u_k + \eta_- u_{k+1} = 0.
\end{equation}
By assumption $u_0 = u_1 = 0$. Hence, applying \eqref{equation:decomp} with $k=1$, we get $\eta_+ u_0 + \phi u_1 + \eta_- u_2 = \eta_- u_2 = 0$ and by injectivity on $\eta_-$ on $\Omega_k$ for $k > 0$ (see Lemma \ref{lemma:injectivity-eta}, case $g=0$), we get $u_2 = 0$. The same argument carries on by induction and shows that $u_k = 0$ for all $k \geq 0$. Similarly, applying \eqref{equation:decomp} for negative $k$ and using the injectivity of $\eta_+$ on $\Omega_k$ for $k < 0$ (see Lemma \ref{lemma:injectivity-eta}, case $g=0$), we get that $u_k = 0$ for all $k < 0$. Hence $u \equiv 0$. Note that this step does not use the positive curvature of the metric.
\end{proof}


 \emph{Step 2. Estimate with compact remainder.} Before proving \eqref{equation:stabilite}, we first prove the following inequality: for all $s \in \R, N > 0$, there exists a constant $C > 0$ such that for all $u \in C^\infty(SM)$,
\begin{equation}
\label{equation:stabilite-rough}
\|u\|_{H^s} \leq C \left( \|(X+\pi^*\phi)u\|_{H^s} +\|u_0\|_{H^{s+1}}+\|u_1\|_{H^{s+1}} + \|u\|_{H^{-N}}\right).
\end{equation}
The proof of \eqref{equation:stabilite-rough} mainly follows that of Theorem \ref{theorem:regularity}. Indeed, we use \eqref{equation:bound3} for $A$ with microsupport and elliptic in a conic neighbourhood of $\HH^*$, and apply the commutation relation \eqref{equation:szego-commutation}:
\begin{equation}
\label{equation:interm}
\begin{split}
\|A \mc{S} u\|_{H^s} &  \leq C \left(\|B_0 (X+\pi^*\phi) \mc{S} u\|_{H^s} + \|\mc{S}  u\|_{H^{-N}} \right) \\
& \leq C \left(\|(X+\pi^*\phi) u\|_{H^s} +\| \eta_- u_0 - \eta_+ u_{-1}\|_{H^s} + \|u\|_{H^{-N}} \right) \\
& \leq C \left(\|(X+\pi^*\phi) u\|_{H^s} +\| u_{-1} \|_{H^{s+1}}+\| u_0\|_{H^{s+1}} + \|u\|_{H^{-N}} \right).
\end{split}
\end{equation}
Observe that
\begin{equation}
\label{equation:decomp2}
\eta_+ u_{-1} + \phi u_0 + \eta_- u_1 = \left( (X+\pi^*\phi) u\right)_0.
\end{equation}
The operator $\eta_+ : C^\infty(M,\Omega_{-1}) \to C^\infty(M,\Omega_0)$ is injective by Lemma \ref{lemma:injectivity-eta}, case $g=0$, and since it is elliptic, it admits a left inverse $L : C^\infty(M,\Omega_0) \to C^\infty(M,\Omega_{-1})$ which is a pseudodifferential operator of order $-1$ such that $L \eta_+ = \mathbbm{1}_{\Omega_{-1}}$ (see Lemma \ref{lemma:ellipticity}, it is possible to remove the smoothing remainder by injectivity of $\eta_+$). Applying $L$ to \eqref{equation:decomp2}, we get:
\[
u_{-1}= L \eta_+ u_{-1}  = - L \phi u_0 - L \eta_- u_1 + L \left( (X+\pi^*\phi) u\right)_0.
\]
Hence, using that $L$ is of order $-1$, we get:
\begin{equation}
\label{equation:estimate-uminus1}
\begin{split}
\|u_{-1}\|_{H^{s+1}} & \leq \|L\phi u_0\|_{H^{s+1}} + \|L \eta_- u_1\|_{H^{s+1}} + \|L \left( (X+\pi^*\phi) u\right)_0\|_{H^{s+1}} \\
& \leq C \left(\|u_0\|_{H^{s+1}} + \|u_{1}\|_{H^{s+1}} + \|(X+\pi^*\phi)u\|_{H^s} \right),
\end{split}
\end{equation}
for some constant $C > 0$, independent of $u$. Reporting this in \eqref{equation:interm} yields:
\[
\|A \mc{S} u \|_{H^s} \leq C \left(\|(X+\pi^*\phi) u\|_{H^s} +\| u_{1} \|_{H^{s+1}}+\| u_0\|_{H^{s+1}} + \|u\|_{H^{-N}} \right)
\]
A similar estimate holds for $A(\mathbbm{1}-\mc{S}) u$, which eventually yields:
\[
\|A u\|_{H^s} \leq C  \left(\|(X+\pi^*\phi) u\|_{H^s}  +  \|u_0\|_{H^{s+1}} +  \|u_{1}\|_{H^{s+1}} + \|u\|_{H^{-N}} \right).
\]

Now, by \S\ref{sssection:projective-positive}, there exists a maximum time $T > 0$ such that for any point $(v,\xi) \in \Sigma$, there exists a time $t \in [-T,T]$ such that $\Phi_t(v,\xi)$ belongs to the elliptic support of $A$. As a consequence, by standard propagation of singularities (see Lemma \ref{lemma:propagation}), there exists $A_{\Sigma} \in \Psi^0(SM)$, microlocally equal to the identity on a conic neighbourhood of the characteristic set $\Sigma$, such that
\begin{equation}
\label{equation:asigma1}
\|A_{\Sigma} u\|_{H^s} \leq C  \left(\|(X+\pi^*\phi) u\|_{H^s}  +  \|u_0\|_{H^{s+1}} +  \|u_{1}\|_{H^{s+1}} + \|u\|_{H^{-N}} \right).
\end{equation}
Since $X+\pi^*\phi$ is elliptic on the microsupport of $\mathbbm{1}-A_{\Sigma}$, by Lemma \ref{lemma:ellipticity}, we can construct $Q \in \Psi^{-1}(SM), R \in \Psi^{-\infty}(SM)$ such that
\begin{equation}
\label{equation:parametrix}
Q(X+\pi^*\phi) = \mathbbm{1}-A_{\Sigma} + R.
\end{equation}
Hence, we get by \eqref{equation:parametrix} that:
\begin{equation}
\label{equation:asigma2}
\begin{split}
\|(\mathbbm{1}-A_{\Sigma})u\|_{H^s} & \leq \|Q(X+\pi^*\phi)u\|_{H^s} + \|Ru\|_{H^s} \\
&  \leq C \left(\|(X+\pi^*\phi)u\|_{H^{s-1}} + \|u\|_{H^{-N}}\right)\\
& \leq C \left(\|(X+\pi^*\phi)u\|_{H^{s}} + \|u\|_{H^{-N}}\right),
\end{split}
\end{equation}
where $C >0$ is a constant independent of $u$ and $N > 0$ can be taken arbitrarily. 
Combining \eqref{equation:asigma1} and \eqref{equation:asigma2} eventually proves \eqref{equation:stabilite-rough}. \\

 \emph{Step 3. Removing the compact remainder.} We now want to get rid of the $\|u\|_{H^{-N}}$ term in \eqref{equation:stabilite-rough}. We argue by contradiction. Assume that \eqref{equation:stabilite} does not hold. Then, for all integer $n > 0$, there exists $u^{(n)} \in H^s(SM)$ such that $\|u^{(n)}\|_{H^s}=1$ and:
\[
\|u^{(n)}\|_{H^s} = 1 \geq n \left( \|(X+\pi^*\phi)u^{(n)}\|_{H^s} + \|u^{(n)}_0\|_{H^{s+1}}+\|u^{(n)}_1\|_{H^{s+1}}\right),
\]
that is, $(X+\pi^*\phi)u^{(n)} \to 0$ in $H^s$ and $u^{(n)}_0,u^{(n)}_1 \to 0$ in $H^{s+1}$. We then apply \eqref{equation:stabilite-rough} with $N \gg 0$ large enough. Note that by compactness of the embedding $H^s(SM) \hookrightarrow H^{-N}(SM)$, we can always assume that, up to extraction, $u^{(n)} \to w$ in $H^{-N}$. As a consequence, we deduce by \eqref{equation:stabilite-rough} that the sequence $(u^{(n)})_{n \geq 0}$ is a Cauchy sequence in $H^{s}(SM)$. By completeness, it implies that $u^{(n)} \to u^\infty$ in $H^s(SM)$ and $\|u^{(n)}\|_{H^s} = 1 = \|u^\infty\|_{H^s}$. Moreover, the distribution $u^\infty$ satisfies $(X+\pi^*\phi)u^\infty = 0$ and $u^\infty_0 = u^\infty_1 = 0$. By Lemma \ref{lemma:inj}, this implies that $u^\infty \equiv 0$, which contradicts $\|u^\infty\|_{H^s}=1$. This eventually proves \eqref{equation:stabilite} and concludes the proof of Theorem \ref{theorem:sphere}.
\end{proof}

%
%

\section{Classification of holomorphic line bundles over twistor space}

\label{section:line-bundles}

\subsection{Holomorphic line bundles over complex manifolds}

\label{section:linebundlesonM}

We first recall some standard results on holomorphic line bundles over a complex manifold $Z$.\footnote{Here $Z$ is a classical complex manifold -- we will explain later how this carries over to twistor spaces, where there is a boundary at which the complex structure degenerates.} Let $L\rightarrow Z$ a complex line bundle.
A {\it partial connection} on $L$ is an operator 
\begin{equation}
\bar \partial_L\colon \Omega^{0,q}(Z,L)\rightarrow \Omega^{0,q+1}(Z,L),\quad q=0,1,
\end{equation}
satisfying the obvious Leibniz rule. The partial connection $\bar \partial_L$ is called {\it integrable}, if $\bar \partial_L^2=0$. 

\begin{definition}
A {\it holomorphic line bundle} on $Z$ is a pair $(L,\bar \partial_L)$ consisting of a complex line bundle $L\rightarrow Z$ and an integrable partial connection 
$\bar \partial_L\colon \Omega^{0,*}(Z,L)\rightarrow \Omega^{0,*+1}(Z,L)$.
Given an open set $U\subset Z$, we say that two holomorphic line bundles $(L,\bar \partial_L)$ and $(L',\bar \partial_{L'})$ are (holomorphically) isomorphic over $U$, written
\begin{equation}
(L,\bar \partial_L) \simeq (L',\bar \partial_{L'}) \text{ on } U,
\end{equation}
if there exists a smooth line bundle isomorphism $\varphi \in C^\infty(U,\mathrm{End}(L,L'))$ such that  $\varphi\bar \partial_{L'}f = \bar \partial_L( \varphi f)$ for all $f\in C^\infty(U,L)$.
\end{definition}

In order to study holomorphic line bundles with a fixed underlying smooth bundle $L\rightarrow Z$, it is convenient to fix an integrable partial connection $\bar \partial_{L}$ as background structure. Every other partial connection is then of the form $\bar \partial_L+\alpha$ for some $\alpha \in \Omega^{0,1}(Z)$ and the integrability condition reads $\bar \partial \alpha = 0 $. Further an isomorphism $(L, \bar \partial_L + \alpha)\simeq (L,\bar \partial_L+\tilde{\alpha})$ on $U\subset Z$ is equivalent to the existence of a function $\varphi \in C^\infty(U,\C^\times)$ with 
\begin{equation}\label{abeq}
\alpha = \varphi^{-1} (\bar \partial + \tilde{\alpha})\varphi.
\end{equation}
Here it is possible to consider $\varphi$ as a $\C$-valued function, because $\mathrm {End}(L,L)$  is canonically isomorphic to $\C$. 

\subsubsection{Picard group} The set of holomorphic line bundles up to global isomorphisms is naturally equipped with a group structure\footnote{The identity element is the trivial line bundle equipped with the trivial flat connection, multiplication is given by the tensor product, and inversion consists in taking the dual.} and is called the \emph{Picard group} $\mathrm{Pic}(Z)$ of $Z$. The subgroup of holomorphic structures on the trivial line bundle $\eps := Z \times \C \to Z$ is denoted by $\mathrm{Pic}_0(Z)$.
A natural background holomorphic structure on $\varepsilon $ is given by the scalar $\bar \partial$-operator. The map
\[
\Omega^{0,1}(Z) \cap \ker \bar \partial \ni \alpha \mapsto L_\alpha := (\eps,\bar \partial + \alpha) \in \mathrm{Pic}_0(Z)
\]
factors through $H_{\bar \partial}^{0,1}(Z)$ since two cohomologous holomorphic $1$-forms yield the same holomorphic structure. Hence,
\begin{equation}
\label{equation:m}
\Phi : H_{\bar \partial}^{0,1}(Z,\C) \ni [\alpha] \mapsto L_\alpha  \in \mathrm{Pic}_0(Z).
\end{equation}
is a well-defined Abelian group homomorphism. This map is surjective but not injective and thus
\begin{equation}
\label{equation:pic}
\mathrm{Pic}_0(Z) \simeq H_{\bar \partial}^{0,1}(Z,\C)/\ker \Phi.
\end{equation} 
The group $\mathrm{Pic}_0(Z)$ can be easily described on all compact Kähler manifolds, hence on all Riemann surfaces. We now assume that $Z := M$ is an oriented Riemann surface of genus $g \geq 0$. The following is standard:

\begin{proposition}
There exists $\omega_1, ..., \omega_{2g} \in H_{\bar \partial}^{0,1}(M,\C)$ independent over the \emph{reals} such that $\ker \Phi$ is given by the lattice $\Z~\omega_1 + ... +\Z~\omega_{2g}$ in $\C^{2g}$. In other words:
\[
\mathrm{Pic}_0(M) \simeq \C^g/\Z^{2g}.
\]
\end{proposition}

The map $\Phi$ is then simply given by the projection map $H_{\bar \partial}^{0,1}(M,\C) \simeq \C^g \to \mathrm{Pic}_0(Z) \simeq \C^g/\Z^{2g}$. For the reader's convenience, we provide a short self-contained proof.

\begin{proof}
Using the Hodge theorem and the fact that $M$ is a complex curve, there is a natural identification
\[
H_{\bar \partial}^{0,1}(M,\C) \simeq \mathscr{H}^{0,1}(M,\C) := \left\{ \omega \in \Omega^{0,1}(M) ~|~ \partial \omega = 0 \right\}.
\]
Under this identification, given $\omega \in H_{\bar \partial}^{0,1}(M,\C)$, one can consider the real-valued form $\eta_\omega := (\omega-\overline{\omega})/i \in H^1_{\mathrm{dR}}(M,\R)$. Note that
\[
\mathscr{I} : H_{\bar \partial}^{0,1}(M,\C) \simeq \C^g \ni \omega \mapsto \eta_\omega \in H^1_{\mathrm{dR}}(M,\R) \simeq \R^{2g}
\]
is only a \emph{real} isomorphism of vector spaces.

Since $d\eta = 0$, the connection $\nabla := d +i \eta$ is a unitary flat connection and thus corresponds to an element of the character variety $\mathrm{Hom}(\pi_1(M),\mathrm{U}(1))$. But the kernel of the map
\[
\Psi : H^1_{\mathrm{dR}}(M,\R) \to \mathrm{Hom}(\pi_1(M),\mathrm{U}(1)),
\]
is easily seen to be given by the integer cohomology $H^1(M,\Z) \simeq \Z^{2g}$ which forms a lattice in $H^1_{\mathrm{dR}}(M,\R)$, that is, it corresponds to all closed forms of the type $\varphi^{-1} d\varphi$, for $\varphi \in C^\infty(M,\Ss^1)$. Moreover, it is immediate to check that the map $\mathscr{I}$ identifies $\ker \Phi$ and $\ker \Psi$, and this proves the claim.
\end{proof}


\subsection{Holomorphic line bundles over twistor space}

We now assume that $Z$ is the twistor space of a closed, oriented Riemannian surface. With Dolbeaut complex defined as in \S\ref{section:dolbeaut}, the definitions from the preceding sections carry over to the present setting.



\begin{proposition}\label{propla1}
Suppose $L_0=(L,\bar \partial_L)$ is a holomorphic line bundle. 
\begin{enumerate}[label=\emph{(\roman*)}]
\item  Let $a\in C^\infty(SM)$ with $a_k =0$ for $k<-1$ and 
\begin{equation}\label{defalpha}
\alpha = \fib_*\left(\sum_{k\ge -1} a_k \omega^{k+1} \tau\right)\in \Omega^{0,1}(Z),
\end{equation}then $L_a = (L,\bar \partial_L + \alpha)$ is a holomorphic line bundle;
\smallskip
\item every holomorphic line bundle on $Z$ that is topologically isomorphic to $L$ is holomorphically isomorphic to $L_a$ for some $a\in \oplus_{k\ge -1}C^\infty(M,\Omega_k)$;
\smallskip
\item we have $L_a \simeq L_0$ on $Z$ if and only if there exists some $\psi \in C^\infty(SM,\C^\times)$ such that $\psi$ and $\psi^{-1}$ are fibrewise holomorphic and
\begin{equation}\label{strongtrivial}
(X+a)\psi=0.
\end{equation}
\end{enumerate}
\end{proposition}

\begin{proof}
	Note that $(\bar \partial_L+\alpha)^2= \bar \partial \alpha$, and $\fib^*\bar \partial \alpha = - \partial_{\bar \omega} \sum_{k\ge -1}a_k\omega^k \cdot \tau\wedge \gamma = 0$ by Proposition \ref{intertwine} and Lemma \ref{smoothtrace}. This shows that $L_a$ is a holomorphic line bundle. For (ii) we may assume that the line bundle in question is of the form $(L,\bar \partial_L + \beta)$ for some $\beta \in \Omega^{0,1}(Z)\cap \ker \bar \partial$.
 Then $\fib^*\beta = h_1 \tau + h_2\gamma$ for $h_1,h_2\in \mathcal I^{-1}(Z)$ and by Lemma \ref{smoothtrace} there exists $g\in \mathcal I^0(Z)$ with $ \partial_{\bar \omega} g = h_2$. Consider
	$
	\varphi = \fib_*(e^{-g}) \in C^\infty(Z,\C^\times),
	$
then
\[
\fib^*\left(\bar \partial\varphi + \beta\varphi\right) = \left[( \omega^2\eta_+ + \eta_-)e^{-g} + h_1e^{-g}\right]\tau + \left[\partial_{\bar \omega} e^{-g} + h_2e^{-g} \right]\gamma
= e^{-g} h\tau, 
\]
where $h = (h_1 - (\omega^2 \eta_+ +\eta_-)g)\in \mathcal I^{-1}(Z)$. We see that \eqref{abeq} is satisfied for $\alpha = \fib_*(h\tau)$ and as $\alpha$ inherits the integrability condition from $\beta$, we must have $0=\fib^*\bar \partial \alpha = - \partial_{\bar \omega} h \cdot \tau \wedge \gamma$. Thus $\partial_{\bar \omega}h=0$ and by Proposition \ref{smoothtrace}, $\alpha$ is of the form \eqref{defalpha} for some $a\in \oplus_{k\ge -1}C^\infty(M,\Omega_k)$.

Finally, if $L$ is topologically trivial, then (iii) is contained in \cite[Proposition 4.11]{Bohr-Paternain-21}. The general case can be treated analogously due to the characterisation \eqref{abeq}.
\end{proof}

If considered up to equivalence on all of $Z$, the moduli space of holomorphic line bundles can be infinite-dimensional. The next proposition makes this precise for the cohomology group $H^{0,1}_{\bar \partial}(Z)$, which should be thought of as tangent space to this moduli space.

\begin{proposition} Let $(M,g)$ be a closed surface for which the tensor tomography problem is solvable\footnote{As mentioned in the introduction, the \emph{tensor tomography problem} consists in studying the equation $Xu = f$, where $f$ has finite Fourier degree and $u \in C^\infty(SM)$. If $f = \pi^* f_0$ for some $f_0\in C^\infty(M)$, then the problem asks whether $f = 0$ and $u$ is constant.} on $C^\infty(M)$ (e.g.~an Anosov surface). Then $H^{0,1}_{\bar \partial}(Z)$ is infinite dimensional.
\end{proposition}

\begin{proof}
We will show that the map $f\mapsto [\fib_*(f\omega\tau)]$ injects $C^\infty(M)\simeq C^\infty(M,\Omega_0)$ into $H^{0,1}_{\bar \partial}(Z)$.
Indeed, if $[\fib_*(f\omega\tau)]=0$, then there exists $h\in \mathcal I^0(Z)\cap \ker \partial_{\bar \omega}$ such that
\begin{equation}
(\omega^2 \eta_+ + \eta_-)h = f\omega,
\end{equation}
where we consider $f$ as function on $SM \times \D$ that is independent of $(v,\omega)$. But then $u=h\vert_{\omega=1}\in C^\infty(SM)$ solves $Xu=f$ on $SM$ and thus $f$ lies in the kernel of the geodesic $X$-ray transform, which is assumed to be trivial.
\end{proof}

\begin{remark}
It is an interesting question whether there are closed surfaces $(M,g)$ for which $H^{0,1}_{\bar \partial}(Z)$ is {\it finite} dimensional. We will not pursue this question here.
\end{remark}
 
\subsection{Isomorphism in the interior} 

A more manageable moduli space of holomorphic line bundles arises by considering isomorphisms only in the interior of $Z$. 
 In this case, it is meaningful (and no more difficult) to consider the twisted bundles $L_a$ from Proposition \ref{propla1} also for certain distributions $a$ and we discuss this first.
 
Let us assume from now on that $(M,g)$ is free of conjugate points such that the unstable bundle $E_u^*\subset T^*(SM)\backslash 0$ is well defined, see \S\ref{prop:coV}. 
Consider a distribution $a\in \mathcal D'(SM)$ with 
\begin{equation}\label{mhodistribution}
\WF(a)\subset E_u^* \quad \text{ and } \quad a_k=0 \text{ for } k<-1.
\end{equation}
As $E_u^* \cap (\HH^* \oplus E_0^*) = \left\{0\right\}$, Proposition \ref{wftest} implies that the Fourier modes of $a$ are smooth, using Proposition \ref{prop_trace1}, one sees that
\begin{equation}
\alpha =\fib_*\left(\sum_{k\ge -1} a_k \omega^{k+1} \tau\right) \in \Omega^{0,1}(Z^\circ)\cap \ker  \bar \partial
\end{equation}
is well defined in the interior of $Z$. Hence, given a background holomorphic line bundle $L_0=(L,\bar \partial_L)$ on $Z$, we obtain a twisted line bundle $L_a=(L,\bar \partial_L + \alpha)$ on $Z^\circ$.


\begin{lemma}\label{lem_wpsi}  
Let $a\in \mathcal D'(SM)$ be as \eqref{mhodistribution}  and suppose there exists a fibrewise holomorphic distribution $w\in \mathcal D'(SM)$ and a function $\psi\in C^\infty(SM,\C^\times)$ such that both $\psi$ and $\psi^{-1}$ are fibrewise holomorphic and such that
\begin{equation}
Xw = a + \psi^{-1} X\psi.
\end{equation}
Then $L_a\simeq L_0$ on $Z^\circ$.
\end{lemma}

\begin{proof}
Let $h=\sum_{k\ge 0}w_k\omega^k\in \mathcal I^0(Z^\circ)$ and $g=\sum_{k\ge 0} \psi_k \omega^k \in \mathcal I^0(Z)$. Then  $\partial_{\bar \omega} h = \partial_{\bar \omega}g=0$ and $g$ is $\C^\times$-valued by the argument principle. We claim that the gauge
\[\varphi = \fib_*\left(e^{-h}g\right)\in C^\infty(Z^\circ,\C^\times)\]
implements the equivalence of $L_a$ and $L_0$ over $Z^\circ$. To this end, compute
\begin{equation*}
\fib^*\left(\varphi^{-1}\bar \partial \varphi\right) = e^h g^{-1} (\omega^2\eta_+ + \eta_-) (e^{-h}g)\cdot \tau  = \left[   - (\omega^2 \eta_+ + \eta_-)h  + g^{-1} (\omega^2\eta_++\eta_-)g\right]\cdot \tau
\end{equation*} 
and
\begin{equation}
(\omega^2\eta_+ + \eta_-)h = \sum_{k\ge 0} (\eta_+w_{k-2} + \eta_-w_k)\omega^k = \sum_{k \ge -1} a_k \omega^{k+1} + \sum_{k\ge -1} (\psi^{-1}X\psi)_k\omega^k.
\end{equation}
Further,
$H=\sum_{k\ge -1} (\psi^{-1}X\psi)_k\omega^k -g^{-1} (\omega^2\eta_++\eta_-)g\in \mathcal I^{-1}(Z)$ is smooth up to the boundary and satisfies
\begin{equation}
\partial_{\bar \omega} H = 0\quad  \text{ and } \quad  H\vert_{\omega=1} =0,
\end{equation}
such that we obtain $H\equiv 0$ from $\mathbf V$-invariance and the identity principle. This shows that
$
\varphi^{-1} \bar \partial \varphi = \alpha
$ and thus $L_a\simeq L_0$ on $Z^\circ$ by \eqref{abeq}.
\end{proof}

\begin{theorem}\label{modulianosov} Suppose $Z$ is the twistor space of an Anosov surface. Let $L_0=(L,\bar \partial_L)$ be a holomorphic line bundle on $Z$ and  let $a\in \mathcal D'(SM)$ be as in \eqref{mhodistribution}. Let $\iota_M\colon M\rightarrow Z$ be the inclusion as $0$-section and consider
\begin{equation}
\alpha = \fib_*\left(\sum_{k\ge -1 } a_k \omega^{k+1} \tau\right) \in \Omega^{0,1}(Z^\circ)\quad \text{ and } \quad \beta = \fib_*\left(\sum_{k\ge 0} a_{2k} \omega^{2k} \bar \tau \wedge \tau\right) \in \Omega^{1,1}(Z^\circ).
\end{equation}
Then $\bar \partial \alpha = 0$ and moreover,
\begin{equation}
\label{equation:eq-th}
\def\arraystretch{1.5}
\left \{ \begin{array}{ll}
\emph{(i)} & [\iota_M^*\alpha]=0 \in  H^{0,1}_{\bar \partial}(M)/H^1(M,\Z)\\
\emph{(ii)} & [\iota_M^*\beta] =0 \in H^2(M)
\end{array} \right.\qquad \Longleftrightarrow \quad L_a\simeq L_0 \text{ on } Z^\circ
\end{equation}
\end{theorem}

\begin{remark}
The implication `$\Longleftarrow$' is easy to prove and holds also in general, if the Anosov property is not satisfied.
\end{remark}

Let us explain the two conditions in more detail. Under the isomorphism between $C^\infty(M,\Omega_{-1})$ $\simeq \Omega^{0,1}(M)$ explained in \textsection \ref{section:fourier}, we have $\iota_M^*\alpha=a_{-1}\in \Omega^{0,1}(M)$ and we may consider its equivalence class $[\iota_M^*\alpha]$ in $H^{0,1}(M)/\ker \Phi$ (with $\ker \Phi \simeq H^1(M,\Z)$) as in \textsection \ref{section:linebundlesonM}. This equivalence class encodes a holomorphic line bundle structure on $M\times \C$ and it is zero, if and only if the holomorphic line bundle is trivial.
This shows that:
\begin{equation}\label{cond1'}
\text{(i)} \quad \Leftrightarrow \quad \exists \psi \in C^\infty(M,\C^\times): ~ (\eta_- + a_{-1})\psi=0.
\end{equation}
Further, by Lemma \ref{lem:volume} we have $\iota_M^*\beta = -2i a_0~ d\mathrm{Vol}_g$, such that the second condition becomes:
\begin{equation}
\text{(ii)} \quad \Leftrightarrow \quad \int_M a_0 =0.
\end{equation}

The proof of the theorem is based on the following consequence of the tensor tomography problem for functions and 1-forms:
\begin{proposition}\label{controlmodes}
 Let $(M,g)$ be an Anosov surface and $q\in C^\infty(M,\Omega_m)$ for $m=0,\pm 1$. Then there exists a distribution $v\in \mathcal D'(SM)$ such that $Xv=0$ and $v_m=q$.
\end{proposition}

\begin{proof} For $m=0$ this is the content of \cite[Theorem 1.4]{Paternain-Salo-Uhlmann-14-2}. For $m=\pm 1$ the claim follows from \cite[Theorem 1.5]{Paternain-Salo-Uhlmann-14-2} or \cite[Corollary 3.8]{Guillarmou-17-1} as we now briefly explain for the $m=1$ case. Take $q\in C^{\infty}(M,\Omega_{1})$ and solve $\eta_{+}p=-\eta_{-}q$ for $p\in C^{\infty}(M,\Omega_{-1})$. This is possible since $\eta_{-}q$ is orthogonal to constants. Then the 1-form represented by $p+q$ is solenoidal and we can apply  \cite[Theorem 1.5]{Paternain-Salo-Uhlmann-14-2} or \cite[Corollary 3.8]{Guillarmou-17-1} to obtain the desired invariant distribution.
\end{proof}

\begin{remark} One can also control two consecutive modes $q_m\in C^\infty(M,\Omega_m)$ and $q_{m+1}\in C^\infty(M,\Omega_{m+1})$ for $m=-1,0$  by first finding $u,v\in \mathcal D'(SM)\cap \ker X$ with $u_m=q_m$ and $v_{m+1}=q_{m+1}$ using the proposition and then putting $w = \dots + v_{m-1}+u_m + v_{m+1} + u_{m+2} +\dots$
\end{remark}

\begin{proof}[Proof of Theorem \ref{modulianosov}] Let us assume that (i) and (ii) hold true. Due to condition (ii) we have $\int_{SM} a = 2\pi \int_M a_0 =0$ and thus the equation $Xu = a$ has a solution $u\in \mathcal D'(SM)$ with $\WF(u)\subset E_u^*$ (this follows from the theory of anisotropic Sobolev spaces for hyperbolic flows, see \cite{Faure-Sjostrand-11, Dyatlov-Zworski-16} for general references and \cite[Theorem 2.6]{Guillarmou-17-1} for a proof of this particular statement). In particular, the Fourier modes of $u$ are smooth and by \eqref{equation:szego-commutation}, we get
\begin{equation}
X\mathcal S u =  a - \eta_+u_{-2} - \eta_+u_{-1}.
\end{equation}
Further, by Proposition \ref{controlmodes} we can find an invariant distribution  $v\in \mathcal D'(SM)$ with $v_{-1} = -u_{-1}$ and $v_{0}=-u_0$. Due to condition (i) in the form \eqref{cond1'} there exists some $\psi\in C^\infty(M,\C^\times)$ with  $ -\psi^{-1} \eta_- \psi = a_{-1} =\eta_+ u _{-2} + \eta_- u_0 $. Hence  
\begin{equation}
X\mc S v = \eta_-v_0  - \eta_+v_{-1} = -\eta_-u_{0} + \eta_+ u_{-1} = \eta_+u_{-2} +  \eta_+ u_{-1}  + \psi^{-1} \eta_- \psi.
\end{equation}	
Setting $w=\mc S u + \mc S v \in \mathcal D'(SM)$, we have shown that
$
X w  = a + \psi^{-1} X\psi
$
and by Lemma \ref{lem_wpsi} this implies $L_a \simeq L_0$ on $Z^\circ$.

For the other direction assume there exists a gauge $\varphi\in C^\infty(Z^\circ,\C^\times)$ such that $-\varphi^{-1}\bar \partial \varphi = \fib_*\left(\sum_{k\ge -1 }a_k\omega^{k+1}\tau\right)$. Then $h = \fib^*\varphi \in \mathcal I^0(Z^\circ)$ satisfies $\partial_{\bar \omega} h = 0$ and thus it can be expanded as a power series $h = \sum_{k\ge 0} h_k \omega^k$ with $h_k\in C^\infty(M,\Omega_k)$. It follows that
\begin{equation}
-\sum_{k\ge 0} (\eta_+h_{k-2} + \eta_-h_k)\omega^k =\left(\sum_{k\ge 0} h_k \omega^k\right)\cdot \left(\sum_{k\ge -1} a_k \omega^{k+1}\right)\quad \text{ on } SM\times \D^\circ,
\end{equation}
where we set $h_k=0$ for $k<0$. Comparing the coefficients of the two power series yields
\begin{equation*}
-\eta_-h_0 =  h_0 a_{-1},\quad -\eta_-h_1 = h_0a_0 + h_1a_{-1},\quad \dots
\end{equation*}
Note that $h_0=\fib^*\varphi\vert_{\omega=0}$ takes values in $\C^\times$, such that \eqref{cond1'} is satisfied with $\psi = h_0$. Further, the preceding display implies that
\begin{equation}
\eta_-(h_0^{-1}h_1) = - h_0^{-2} (\eta_-h_0)h_1 + h_0^{-1}\eta_- h_1 = h_0^{-1}(\eta_-h_1 + a_{-1}h_1) = -a_0.
\end{equation}
Recall that $\eta_-=\frac 12(X+iH)$ and that both $X$ and $H$ preserve the Liouville measure, such that $2\pi \int_M a_0 = -\int_{SM} \eta_-(h_0^{-1}h_1) = 0$.
\end{proof}

\appendix

\section{A hitchhiker's guide to microlocal analysis}

\label{appendix}

We refer to \cite{Hormander-90, Hormander-94, Shubin-01, Grigis-Sjostrand-94,Lefeuvre-book} for an extensive treatment of pseudodifferential operators on closed manifolds.

\subsection{Definitions. First properties. Ellipticity.}

Let $\M$ be a closed  manifold of dimension $n$. Define $S^k(T^*\M) \subset C^\infty(T^*\M)$, for $k \in \R$, the space of symbols as the set of smooth functions $a$ on $T^*M$ satisfying the following bounds, in any coordinate chart $U \subset \R^n$: for all compact subset $K \subset U$, $\alpha,\beta \in \N^n$, there exists a constant $C > 0$ such that
\begin{equation}
\label{equation:symbol}
\forall (x,\xi) \in T^* K \simeq K \times \R^n, \qquad |\partial^\alpha_\xi \partial^\beta_x a(x,\xi)| \leq C \langle \xi \rangle^{k-|\alpha|}.
\end{equation}
Since \eqref{equation:symbol} is invariant by diffeomorphism on $\R^n$, $S^k(T^*\M)$ is intrinsically defined on the manifold $\M$.

The set of \emph{smoothing operators} $\Psi^{-\infty}(\M)$ is the space of linear operators on $\M$ with smooth Schwartz kernel. We denote by $\Op$ a quantization on $\M$, given in a local coordinate chart $U \subset \R^n$ by: 
\[
\Op(a)f (x) = \dfrac{1}{(2\pi)^n} \int_{\R^n_\xi} \int_{\R^n_y} e^{i\xi\cdot(x-y)} a(x,\xi) f(y) |dy| |d\xi|,
\]
where $a \in S^k(T^*U), f \in C^\infty_{\comp}(U)$. We define \emph{pseudodifferential operators} of order $k \in \R$ as the set of operators
\[
\Psi^k(\M) := \left\{ \Op(a) + R ~|~ a \in S^k(T^*\M), R \in \Psi^{-\infty}(\M)\right\}.
\]
It is standard that $\Psi^k(\M)$ is intrinsically defined, independently of the choice of quantization $\Op$ (see \cite[Chapter 18]{Hormander-94} or \cite[Chapter 5]{Lefeuvre-book} for the theory of pseudodifferential operators and their symbols).

The \emph{principal symbol map}
\[
\sigma : \Psi^k(\M) \to S^k(T^*\M)/S^{k-1}(T^*\M)
\]
is defined so that $\Op(\sigma(A))-A \in \Psi^{k-1}(\M)$ for all $A \in \Psi^k(\M)$, that is, we have the exact sequence:
\[
0 \longrightarrow \Psi^{k-1}(\M) \longrightarrow \Psi^k(\M) \longrightarrow S^k(T^*\M)/S^{k-1}(T^*\M) \longrightarrow 0.
\]
For $A \in \Psi^k(T^*\M)$, its elliptic set $\mathrm{ell}(A) \subset T^*\M \setminus \left\{0\right\}$ is defined as the open conic set of points $(x_0,\xi_0) \in T^*\M \setminus \left\{ 0\right\}$ such that there exists $\eps > 0$ with the following property:
\begin{equation}
\label{equation:elliptic}
\left(|\xi| \geq 1/\eps \text{ and } d_{S^*\M}\left((x,\xi/|\xi|), (x_0,\xi_0/|\xi_0|) \right) < \eps \right) \implies |\sigma_A(x,\xi)| \geq \eps \langle \xi\rangle^k.
\end{equation}
Here $d_{S^*\M}$ is any metric on the cosphere bundle $S^*\M := T^*\M/\R_+$, where the $\R_+$-action is given by radial dilation in the fibres of $T^*\M$. An operator is said to be \emph{elliptic} if $\mathrm{ell}(A) = T^*\M \setminus \left\{0\right\}$. The \emph{characteristic set} $\Sigma(A)$ of an operator is the closed conic subset defined as the complement of the elliptic set in $T^*\M$.

\begin{example}
\label{example:vfield}
Let $X \in C^\infty(\M,T\M)$ be a vector field on $\M$, seen as a differential operator of order $1$. Then, its principal symbol is $\sigma_X(x,\xi) = i\langle \xi,X(x)\rangle$ and thus, using \eqref{equation:elliptic}, it is immediate to check that $\Sigma(X) = \left\{\langle\xi,X(x)\rangle=0\right\}$ and $\mathrm{ell}(X) = \left(T^*\M \setminus \left\{0\right\}\right) \setminus \Sigma(X)$.
\end{example}

The \emph{wavefront set} $\WF(A)$ (or the \emph{microsupport}) of an operator $A \in \Psi^{k}(\M)$ is the (closed) conic subset of $T^*\M \setminus \left\{0\right\}$ satisfying the following property: $(x_0,\xi_0) \notin \WF(A)$ if and only if for all $m \in \R$, for all $b \in S^m(T^*\M)$ supported in a small conic neighbourhood of $(x_0,\xi_0)$, one has $A\Op(b) \in \Psi^{-\infty}(\M)$. In other words, the complement of the wavefront set of $A$ is the set of codirections where $A$ behaves as a smoothing operator.

\begin{example}
Let $\M = \R^n$ and $a(\xi) := \langle\xi\rangle^m\chi(\xi)$, where $\chi \in C^\infty(\R^n)$ is a smooth function such that $\chi(\lambda \xi) = \chi(\xi)$ for all $\lambda \geq 1, |\xi| \geq 1$ (that is $\chi$ is $0$-homogeneous outside of the unit ball). Observe that $a \in S^m(T^*\R^n)$. Let $\mc{C}_\chi$ denote the cone generated by $\supp(\chi|_{S^n})$, that is $\mc{C}_\chi :=\{\lambda \xi ~|~ \lambda > 0, \xi \in \mathrm{supp}(\chi), |\xi|=1\}$. Define
\[
A f(x) := \Op(a)f(x) = \mc{F}^{-1}(a\mc{F} f)(x) = \dfrac{1}{(2\pi)^n} \int_{\R^n \times \R^n} e^{i\xi(x-y)} a(\xi) f(y) \dd y \dd \xi,
\]
where $\mc{F}$ denotes the Fourier transform. We claim that
\begin{equation}
\label{equation:wfa}
\WF(A) = \R^n \times \mc{C}_\chi \subset T^*\R^n \simeq \R^n \times \R^n.
\end{equation}
Indeed, since $A$ commutes with all translations in $\R^n$, its wavefront set must be translation invariant in the $x$-variable hence of the form $\WF(A) = \R^n \times \mc{B}$ for some cone $\mc{B} \subset \R^n$. Let $b \in S^0(T^*\R^n)$ be a symbol depending only on the $\xi$-variable, $0$-homogeneous for $|\xi| \geq 1$ (similarly to $\chi$), such that $\mc{C}_b \cap \mc{C}_\chi = \emptyset$. Then, by construction, $A\Op(b) = \Op(\langle\xi\rangle^m\chi b) \in \Psi^{-\infty}(\R^n)$ since $\chi b \equiv 0$ outside of the unit ball. This shows that $\WF(A) \subset \R^n \times \mc{C}_\chi$. To prove the converse inclusion, one can exhibit a distribution $u \in \mc{D}'(\R^n)$ such that $\WF(u) \subset \R^n \times \mc{C}_\chi$ (e.g. a Dirac mass on a hyperplane in $\R^n$ whose normal vector is contained in $\mc{C}_\chi$) and $A u \notin C^\infty(\R^n)$. This is an equivalent characterization of the wavefront set of an operator and finally proves \eqref{equation:wfa}.
\end{example}

Pseudodifferential operators in $\R^n$ obtained by quantization of symbols depending only on the $\xi$-variable (as in the previous example) are called \emph{Fourier multipliers}. It can be checked that they coincide exactly with the set of pseudodifferential operators commuting with all translations in $\R^n$.

The wavefront set $\WF(u)$ of a distribution $u \in \mc{D}'(SM)$ is the (closed) conic subset of $T^*\M \setminus \left\{0\right\}$ satisfying the following property: $(x_0,\xi_0) \notin \WF(u)$ if and only if there exists a small open conic neighbourhood $V$ of $(x_0,\xi_0)$ such that for all $k \in \R$, for all $A \in \Psi^k(\M)$ with wavefront set contained in $V$, one has $A u \in C^\infty(\M)$.

\begin{example}
Let $(x\pm i0)^{-1}$ denote the distribution defined for $\varphi \in C^\infty_{\comp}(\R)$ by
\[
((x\pm i0)^{-1},\varphi) := \lim_{\eps \to 0} \int_{\R} \varphi(x) (x\pm i\eps)^{-1} \dd x.
\]
We claim that $\WF((x \pm i0)^{-1}) = \{0\} \times (0,\pm \infty) \subset T^*\R \simeq \R \times \R$. Indeed, $(x\pm i0)^{-1}$ coincides with the smooth function $\R \setminus \{0\} \ni x \mapsto 1/x$ outside of the origin so its wavefront set is contained in $\{0\} \times \R$. Additionally, it is a distribution of order $1$, so it does not coincide with any smooth function near the origin, and its wavefront set is therefore non empty.
 
Writing
\begin{equation}
\label{equation:ec}
(x+i0)^{-1} = \lim_{\eps \to 0} \dfrac{1}{i} \int_0^{+\infty} e^{-\eps \tau} e^{i x \tau} \dd\tau,
\end{equation}
where the limit is taken in $\mc{D}'(\R)$ and using the definition of the wavefront set, it can be verified that $\{0\} \times (0,\mp \infty) \notin \WF((x \pm i0)^{-1})$, and this proves the claim. (This is because \eqref{equation:ec} shows that $(x+i0)^{-1}$ is a sum of positive Fourier modes, so applying a Fourier multiplier $\Op(a)$ with $a \in C^\infty(\R_-)$, $a(\xi) \equiv 1$ for all $\xi \leq -1$ and $a(\xi) \equiv 0$ on $[-1/2,0]$, one finds that $\Op(a)(x+i0)^{-1} \in C^\infty(\R)$.)

Let us provide a more conceptual proof of this fact in the spirit of the present paper. The biholomorphism $\Phi : \D \to \HH = \{z \in \C ~|~ \Im(z) > 0\}$ defined by $\Phi(y) := i(1+y)(1-y)^{-1}$ identifies the open unit disk with the Poincaré half-space. The function $\Phi^*(1/z)$ is a holomorphic function in $\D$ blowing up polynomially near the boundary $\partial \D \simeq S^1$ as in \eqref{equation:blowup-pol} (more precisely, it only blows up near the point $-1 = \Phi^{-1}(0)$). Hence, the trace of $\Phi^*(1/z)$ on $S^1$ defines a distribution $u \in \mc{D}'(S^1)$ and it is straightforward to check that $\Phi^*((x+i0)^{-1}) = u$. Following the arguments of \eqref{equation:wf-szego2}, one obtains that $\WF(u) \subset \{ (v,\xi) ~|~ v \in S^1, \xi \in T^*_vS^1, (\xi,\partial_\theta) > 0\}$, where $\partial_\theta$ generates the rotation on $S^1$. Finally, using that wavefront set is preserved under the action of diffeomorphisms, one finds that $\WF((x+i0)^{-1}) = \{0\} \times (0,\infty)$.
\end{example}

We refer to \cite[Chapter 3]{Lefeuvre-book} for an extensive discussion on the notion of wavefront set. \\

The important property of elliptic operators is that they are invertible modulo smoothing remainders:

\begin{lemma}
\label{lemma:ellipticity}
Let $A \in \Psi^k(\M)$. Then, for all $T \in \Psi^0(T^*\M)$ such that $\WF(T) \subset \mathrm{ell}(A)$, there exists $B \in \Psi^{-k}(T^*\M), R_L, R_R \in \Psi^{-\infty}(\M)$ such that
\[
B A = T+ R_L, \qquad AB = T + R_R.
\]
In particular, if $\mathrm{ell}(A) = T^*\M \setminus\left\{0\right\}$, one can take $T=\mathbbm{1}$.
\end{lemma}

We refer to \cite[Section 4.2.1.2]{Lefeuvre-book} for a proof of the previous lemma.

\subsection{Sobolev spaces. Elliptic estimate} Let $g$ be an arbitrary metric on $\M$ and let $\Delta_g \geq 0$ be the nonnegative Hodge Laplacian acting on functions. For $s \in \R$, define $(\mathbbm{1}+\Delta)^s$ by the spectral theorem. This is an invertible pseudodifferential operator of order $2s$.

For $s \in \R$, $u \in C^\infty(\M)$, we set
\begin{equation}
\label{equation:norm}
\|u\|_{H^s} := \|(\mathbbm{1}+\Delta)^{s/2}u\|_{L^2},
\end{equation}
and define $H^s(\M)$ to be the completion of $C^\infty(\M)$ with respect to \eqref{equation:norm}. The space $H^s(\M)$ is intrinsically defined, independently of the choice of metric $g$, and changing the metric only changes the norm \eqref{equation:norm} by an equivalent norm. A distribution $u \in \mc{D}'(SM)$ is said to \emph{microlocally $H^s$ near $(x_0,\xi_0) \in T^*\M \setminus \left\{0\right\}$} if the following holds: there exists a small open conic neighbourhood $V$ of $(x_0,\xi_0)$ such that for all $A \in \Psi^0(\M)$ with $\WF(A) \subset V$, one has $A u \in H^s(\M)$.

The following boundedness result for pseudodifferential operators holds: for all $k \in \R, A \in \Psi^{k}(\M)$ and $s \in \R$,
\begin{equation}
\label{equation:boundedness}
A : H^{s+k}(\M) \to H^s(\M)
\end{equation}
is bounded. In particular, combining the parametrix construction of Lemma \ref{lemma:ellipticity} for elliptic operators and the boundedness \eqref{equation:boundedness}, one obtains the following:

\begin{lemma}
Let $A \in \Psi^k(\M)$ be an elliptic pseudodifferential operator. Then, for all $s \in \R, N > 0$, there exists a constant $C >0$ such that: for all $u \in C^\infty(\M)$,
\begin{equation}
\label{equation:elliptic-bound}
\|u\|_{H^{s+k}} \leq C \left( \|Au\|_{H^s} + \|u\|_{H^{-N}}\right).
\end{equation}
Moreover, if $u \in \mc{D}'(\M)$ is merely a distribution, $u \in H^{-N}(\M)$ and $Au \in H^s(\M)$, then $u \in H^{s+k}(\M)$ and \eqref{equation:elliptic-bound} holds.
\end{lemma}

We now describe a similar bound to \eqref{equation:elliptic-bound} in the case where $A$ is not elliptic but of \emph{real principal type}. This is known as the \emph{propagation of singularities} for pseudodifferential operators.

\subsection{Propagation of singularities} Let $P \in \Psi^m(\M)$ be a pseudodifferential operator. We will say that $P$ is of \emph{real principal type} if its principal symbol is real-valued and homogeneous (of order $m$). For such an operator, we denote by $H_P \in C^\infty(T^*\M,T(T^*\M))$ the Hamiltonian vector field on $T^*\M$ (equipped with the standard Liouville $2$-form) generated by the principal symbol $\sigma_P$, and by $(\Phi_t)_{t \in \R}$ the Hamiltonian flow it generates.

\begin{lemma}
\label{lemma:propagation}
Let $P \in \Psi^m(\M)$ be a pseudodifferential operator of real principal type. Let $A, B, B_0 \in \Psi^0(\M)$ such that the following holds: for every $(x,\xi) \in \WF(A)$, there exists a time $T > 0$ such that $\Phi_{-T}(x,\xi) \in \mathrm{ell}(B)$ and for all $t \in [0,T]$, $\Phi_{-t}(x,\xi) \in \mathrm{ell}(B_0)$. Then, for all $s \in \R, N > 0$, there exists a constant $C > 0$ such that: for all $u \in C^\infty(\M)$,
\begin{equation}
\label{equation:propagation}
\|Au\|_{H^{s+m-1}} \leq C \left( \|Bu\|_{H^{s+m-1}} + \|B_0 P u\|_{H^s} + \|u\|_{H^{-N}} \right).
\end{equation}
Moreover, if $u \in \mc{D}'(\M)$ is merely a distribution, $u \in H^{-N}(\M)$ , $Bu \in H^{s+m-1}(\M)$ and $B_0P u \in H^s(\M)$, then $u \in H^{s+m-1}(\M)$ and \eqref{equation:propagation} holds.
\end{lemma}

We refer to \cite[Theorem E.47]{Dyatlov-Zworski-book} for a proof. In particular, if $Pu = 0$, then Lemma \ref{lemma:propagation} implies that $\WF(u)$ is invariant under the flow $(\Phi_t)_{t \in \R}$ which is the content of \cite[Theorem 26.1.1]{Hormander-4}. Lemma \ref{lemma:propagation} is applied in this article with $P := \pm iX \in \Psi^1(\M)$, where $X$ is a smooth vector field on $\M$. By Example \ref{example:vfield}, $\sigma_P(x,\xi) = \langle \xi,X(x) \rangle$ and the Hamiltonian flow generated by $H_P$ is simply given by $\Phi_t(x,\xi) = (\varphi_t(x), d\varphi_t^{-\top}(x)\xi)$, the symplectic lift of the flow $(\varphi_t)_{t \in \R}$ generated by $X$ to $T^*\M$. (Here, ${}^{-\top}$ stands for the inverse transpose.)

\bibliographystyle{plain}
\bibliography{JLMS_revision.bbl}

\end{document}